\documentclass[12]{article}
\usepackage{amssymb,amsmath,amsthm,amsfonts,amssymb}
\usepackage{hyperref}

\title{\textsf
{Generic complexity of the Conjugacy Problem in HNN-extensions and
algorithmic stratification of Miller's groups}}
\date{Version 1, December 1,   2006}
\author{\textsf{Alexandre V. Borovik}\thanks{Partially supported
by the Royal Society Leverhulme Trust Senior Research Fellowship.}
\and \textsf{Alexei G. Myasnikov}\thanks{Supported by NSF grant
DMS-0405105,
 NSERC Discovery grant RGPIN 261898, and NSERC Canada Research Chair grant.} \and \textsf{Vladimir N.
Remeslennikov\thanks{Supported by EPSRC grant GR/R29451 and by RFFI
grant 02-01-00192.}}}

\date{}

\newcommand{\bea}{\begin{eqnarray*}}
\newcommand{\eea}{\end{eqnarray*}}
\newcommand{\bq}{\begin{quote}}
\newcommand{\eq}{\end{quote}}
\newcommand{\beq}{\begin{equation}}
\newcommand{\eeq}{\end{equation}}
\newcommand{\bi}{\begin{itemize}}
\newcommand{\ei}{\end{itemize}}

\newcommand{\bd}{\begin{description}}
\newcommand{\ed}{\end{description}}

\newcommand{\BH}{\ensuremath{\mathbb B\mathbb H}}
\newcommand{\SBH}{\ensuremath{\mathbb S\mathbb B\mathbb H}}

\newtheorem{corollary}{Corollary}[section]
\newtheorem{theorem}[corollary]{Theorem}
\newtheorem{lemma}[corollary]{Lemma}

\newtheorem{proposition}[corollary]{Proposition}

\theoremstyle{definition}

\newtheorem{example}[corollary]{Example}

\begin{document}

\maketitle

\begin{center}
{\it To Boris Plotkin as a sign of our friendship and respect.}
\end{center}

\begin{abstract} We discuss  time complexity of The Conjugacy Problem in
HNN-extensions of groups, in particular, in Miller's groups. We show
that for ``almost all'', in some explicit sense, elements, the
Conjugacy Problem is decidable in cubic time. It is worth noting
that the Conjugacy Problem in  a Miller group may have be
undecidable. Our results show that ``hard'' instances of the problem
comprise a negligibly small part of the group.
\end{abstract}

\tableofcontents

\section{Introduction}

The present paper is concerned with the generic complexity of the
Conjugacy Problem in HNN-extensions of groups, in particular, in
Miller's groups. Starting with a presentation for a finitely
presented group $H$, Miller \cite{miller2} constructed a generalized
HNN-extension $G(H)$ of a free group; he then showed that the
Conjugacy Problem in $G(H)$ is undecidable provided  the Word
Problem is undecidable in $H$. Varying the group $H$, one can easily
construct infinitely many groups $G(H)$ with decidable word problem
and undecidable conjugacy problem. Moreover, even the class of free
products $A \ast_C B$ of free groups $A$ and $B$ with amalgamation
over a finitely generated subgroup $C$ contains specimens with
algorithmically undecidable conjugacy problem \cite{miller1}.

This remarkable result shows that the conjugacy problem can be
surprisingly difficult even in groups whose structure we seem to
understand well. In the next few years more examples of
HNN-extensions with decidable word problem and undecidable conjugacy
problem followed (see, for example,  \cite{bk}). Striking
undecidability results of this sort scared away any general research
on the word and conjugacy problems in amalgamated free products and
HNN-extensions. The classical tools of amalgamated products and
HNN-extensions have been abandoned and replaced by methods of
hyperbolic groups \cite{bf,km1,Mikhailovskii}, or automatic groups
\cite{bgss,Eps}, or relatively hyperbolic groups \cite{Bum,Osin}.

In this and other papers in a series of works on algorithmic
problems in amalgamated free products and HNN-extensions
\cite{amalgam-1,complexity,HNN1} we make an attempt to rehabilitate
the classical algorithmic techniques to deal with amalgams.  Our
approach treats  both decidable and undecidable cases
simultaneously, as well as the case of hyperbolic groups mentioned
above. We show that, despite the common belief, the Word and
Conjugacy Problems in amalgamated free products and HNN-extensions
of groups are generically easy
  and the classical algorithms are very fast  on ``most"
 or ``typical" inputs.  In fact, we analyze the computational complexity
 of even harder  algorithmic problems which
lately attracted much attention in cryptography (see
\cite{AAG,KL,petrides}, and surveys \cite{DHorn,S}), the so-called
{\em Normal Form Search Problem} and {\em Conjugacy Search Problem}.
The former one requires for a given element $g$ of a group $G$ to
find the unique normal form of $G$ (assuming that the normal forms
of elements of $G$ are fixed in advance). The latter asks for an
algorithm to check whether or not two given elements of $G$ are
conjugate in $G$, and if they are, to find a conjugator. Our
analysis is based on recent ideas of stratification and generic
complexity \cite{multiplicative,KMSS}; the appendix to the paper
contains the necessary definitions from \cite{multiplicative} on
asymptotic classification of subsets in groups.

Although  the present paper is essentially independent from the
other papers in the series \cite{amalgam-1,complexity,HNN1}, it
might be useful to discuss some of their results.

In \cite{amalgam-1,complexity}, working under some mild assumptions
about the groups involved in a given free amalgamated product of
groups $G$, we stratify $G$ into two parts with respect to the
``hardness'' of the conjugacy problem:
 \begin{itemize}
\item  a \emph{Regular Part} $RP$, consisting  of so-called {\em
regular elements} for which the conjugacy problem is decidable by
 standard algorithms. We show that the regular part $RP$
has very good algorithmic properties:

 \begin{itemize}
 \item the standard algorithms are very fast on regular elements;

 \item  if an element is a conjugate of a given
 regular element then the algorithms quickly provide a conjugator,
 so the Conjugacy Search problem is also decidable for
regular elements;
 \item the set $RP$ is \emph{generic} in $G$,
that is, it is very ``big'' (asymptotically the whole group, see
Sections  \ref{subesc:measuring} and \ref{subsec:negligible});
 \item RP is decidable;
\end{itemize}
 \item  the {\em Black Hole} $BH$ (the complement of the set of regular elements)
  which consists of elements in $G$ for
 which either the standard algorithms  do not work at all, or they
 require a considerable modification, or it is not clear yet whether
 these algorithms work or not.

  \end{itemize}

In this paper we show that similar results hold for HNN-extensions
of groups.  This general technique for solving the conjugacy problem
in HNN-extensions does not work in those, very rare, groups where
the Black Hole ($BH$) of the conjugacy problem coincides with the
whole group, in particular in Miller's groups (see Lemma
\ref{le:BH}). However, the conjugacy problem in Miller's groups is
still easy for most of the elements in $BH$. In this case one has to
stratify the Black Hole itself. To this end, we introduce the notion
of a Strongly Black Hole $\SBH$ (see Section \ref{sec:6}). It is
proven that the Conjugacy Search Problem for elements that do not
lie in the Strong Black Hole $\SBH$ is decidable in cubic time
(Theorem \ref{th:main}). We give an explicit description of the size
of $\SBH$ for Miller's groups and prove that $\SBH$ is a strongly
negligible set (Theorem \ref{thm:8.1}).

This is the first example of a non-trivial solution of the
Stratified Conjugacy Problem in a finitely presented group with
undecidable conjugacy problem.

Throughout the paper we mention various algorithmic problems in
groups. A suitable  discussion on this can be found in
\cite{amalgam-1}.

\section{HNN-extensions}

\subsection{Preliminaries}

We introduce in brief some terminology and formulate several known
results on HNN-extensions of groups. We refer to the books
\cite{LS,miller1} and one of the original papers \cite{collins}
for more detail.

 Let  $H= \left< X \mid
\mathcal{R}\right>$ be a group given by
 generators and relators, and let $A = \left< U_i \ \mid \  i \in I\right>$ and
 $B = \left< V_i \ \mid \  i \in
I\right>$ be two isomorphic subgroups of $H$ generated,
correspondingly,  by elements $U_i$  and $V_i$ $(i \in I)$ from $H$
which are given as words in $X \cup X^{-1}$. Let
$$\phi: A \rightarrow B$$
 be an isomorphism  defined  by $\phi:U_i \rightarrow V_i$, $i
\in I$.  Then the group  $G$ defined by the presentation
\[
G = \left< X, t \mid \mathcal{R},\, t^{-1}U_it = V_i,\, i \in I
\right>
\]
is called an {\em HNN-extension} of the {\em base} group $H$ with
the {\em stable} letter $t$ and {\em associated} (via the
isomorphism $\phi$) subgroups $A$ and $B$. We sometimes write $G$
as
\[
G = \left< H, t \mid t^{-1}At = B, \, \phi\right>.
\]
An HNN-extension $G$ is called \emph{degenerate} if $H=A=B$.

 A
modification of the above definition is that of \emph{multiple
HNN-extension}. The data consist of a group $H$ and a set of
isomorphisms $\phi_i: A_i \rightarrow B_i$ between subgroups of $H$.
Then extending the case above we define a multiple HNN-extension of
$H$ as
 $$G = \left< H, t_i \mid t_i^{-1}A_it_i = B_i,
\phi_i,\, (i\in I)\right>.$$

\subsection{Reduced and normal forms}
\label{sec:reducedandnormal}

The main focus of this section  is on algorithms for computing
\emph{reduced}
 and \emph{normal forms} of elements in HNN-extensions of groups. We consider only
HNN-extensions with one stable letter, but one can easily extend
the results to arbitrary multiple HNN-extensions.

Let $ G = \left< H, t \mid t^{-1} At = B, \, \phi\right>$  be an
HNN-extension of a group $H$ with  stable letter $t$ and associated
subgroups $A, B$. Every element $g$ of $G $ can be written in the
form
 \begin{equation}
 \label{eq:form}
 g = w_0t^{\epsilon_1}w_1 \cdots
t^{\epsilon_n}w_n,
 \end{equation}
 where $\epsilon_i = \pm 1$ and
$w_i$ is a (possibly empty) word in the generating set $X$. The
following result is well known (see, for example, \cite{LS}).
\begin{theorem} 
\label{th:4}
 Let $G=\left< H,t \mid t^{-1}At= B, \phi \right>$,
and let $$g = w_0t^{\epsilon_1}w_1 \cdots t^{\epsilon_n}w_n.$$
If\/ $g$ represents the identity element of $G$ then either
\begin{itemize}
\item[{\rm (a) }]  $n= 0 $ and $w_0$ represents the identity
element of $H$; or

\item[{\rm (b) }]  $g$ contains a subword of the form either
$t^{-1}w_it$ with $w_i \in A$ or $tw_it^{-1}$ with $w_i \in B$
{\rm (}words of this type are called\/ \emph{pinches}{\rm )}.
\end{itemize}
\end{theorem}

Theorem \ref{th:4} immediately gives a decision algorithm for the
Word Problem in $G$ provided one can  effectively solve   the Word
Problem  ``Is $w_0 =1$?'' and Membership Problems  ``Are $w_i
 \in A$ and/or $w_i \in B$?'' in the
group $H$.  We will
 have to say more on the time complexity of the Word Problem in
 $G$ in the sequel.

 We say that (\ref{eq:form}) is  a \emph{reduced form}  of $g
\in G$ if no pinches occur in it. It can be shown that the number
of occurrences of $t_i$ in a reduced form of $g$ does not depend
on the choice of reduced form; we shall call it the \emph{length}
of $g$ and denote it by $l(g)$.

 We say that an element $g$ with $l(g)
>0 $ is \emph{cyclically reduced} if $l(g^2) = 2l(g)$. In addition,
we impose extra conditions in case
$l(g) =0 $ (which is equivalent to saying that $g \in H$): namely,
we say that $g$ is cyclically reduced if either $g \in A \cup B$
or $g$ is not conjugate in $H$ to any element from $A \cup B$.

Equivalently, the definition of cyclically reduced elements can be
formulated as follows. A reduced form
\[
g = h t^{\epsilon_1}s_1 \cdots t^{\epsilon_n}s_n
\]
of an element $g$ is {\em cyclically reduced} if and  only if
 \begin{itemize}
 \item If $n=0$ then either $h\in A\cup B$ or $h$ is not conjugate
in $G$ to any element in $A \cup B$.
 \item if $n >0$ then either $\epsilon_1 = \epsilon_n$, or $s_nh$
does not belong to $A$ provided  ${\epsilon_n} = -1$, or  $s_nh$
does not belong to $B$ provided  ${\epsilon_n} = 1$.

\end{itemize}

We warn that our definition of cyclically reduced elements differs
from that of \cite{LS}; elements cyclically reduced in our sense are
 cyclically reduced in the sense of \cite{LS} but not vice-versa.

 Cyclically reduced forms of elements in $G$ are not unique. To define unique
  \emph{normal forms} of elements in $G$ one needs to fix
  systems of right coset representatives of $A$ and $B$ in $G$.

Let $S_A$ and $S_B$ be systems of right coset representatives
(transversals) of the subgroups $A$ and $B$ in $H$ (we always assume
that the identity element $1$ is the representative of $A$ and $B$
in $H$). A reduced form
  \begin{equation}
  \label{eq:normalforms}
  g = h_0t^{\epsilon_1}h_1 \cdots t^{\epsilon_n}h_n
  \end{equation}
of an element $g \in G$ is said to be a  \emph{normal form} of $g$
if the following conditions hold:
 \begin{itemize}
  \item  $h_0\in H$;
   \item  if $\epsilon_i = -1$ then $h_i \in S_A$;
   \item  if $\epsilon_i = 1$ then $h_i \in S_B$.
\end{itemize}
Normal forms of elements of $G$ are unique  in the sense that  the
elements $h_0, \ldots,
 h_n \in H$ in (\ref{eq:normalforms})  are uniquely defined by $g$ (see, for
 example, \cite{LS}). However, they could be presented by different
 words in the generating set $X$  of $H$. To require uniqueness of
 representation (\ref{eq:normalforms}) by words in $X$
  one has to assume that elements of $H$ can be uniquely presented
  by some particular words in $X$,  i.e.,  existence of normal forms of elements in $H$.

It is convenient sometimes to write down
 the normal form (\ref{eq:normalforms}) of $g$ as
  \begin{equation} \label{eq:normal-2}
g =  h_0p_1\cdots p_k
\end{equation}
where  $p_i = t^{\epsilon_i}s_i$ and $ s_i \in S_A$ if $
\epsilon_i = -1$, $ s_i \in S_B$  if $ \epsilon_i = 1.$ Observe
that this decomposition corresponds to the standard decomposition
of elements of $G$ when $G$ is viewed as  the universal Stallings
group $U(P)$ associated with the pregroup
\[P = \{H, tH, t^{-1}H\},\] (for a more detailed
description of pregroups see \cite{MRS}).

\subsection{Algorithm I for computing reduced forms}
This algorithm takes as  input a word of the form
\[
g = w_0t^{\epsilon_1}w_1 \cdots t^{\epsilon_n}w_n.
\]
If the word contains no pinches then it is reduced. Otherwise, we
look at the first on the left subword of the form
$t^{\epsilon_i}w_it^{\epsilon_{i+1}}$ that is a pinch and transform
the subword according to one of the rules:

\bi \item If $w_i \in A$ and $\epsilon_i = -1$ then rewrite $w_i$ in
the given generators $U_j, j \in I,$ for $A$ and replace
$t^{-1}w_it$ by $\phi(w_i)$, using substitution $t^{-1}U_jt
\rightarrow V_j$;

\item If $w_i \in B$ and $\epsilon_i = 1$ then rewrite $w_i$ in
the given generators $V_j, j \in I,$ for $B$ and  replace
$tw_it^{-1}$ by $\phi^{-1}(w_i)$, using substitution $tV_jt^{-1}
\rightarrow U_j$; \ei
 thus decreasing the  length
$l(g)$ of the word by $2$. Notice that to carry out this algorithm
one needs to be able to verify whether or not an element $w \in H$,
given as a word in the generators of $H$,  belongs to the subgroup
$A$ or $B$, and, if it does, then  to rewrite $w$ as a word in the
given generators of $A$ or $B$. In this event we say that  the {\em
Search Membership Problem} ($\mathbf{SMP}$) is  decidable for the
subgroups $A$ and $B$ in $H$.

We summarize this discussion in the following result (similar to the
one for amalgamated products \cite{amalgam-1}).

\begin{proposition}
Let $G=\langle H,t \mid t^{-1} A t=B\rangle$ be an HNN-extension of
a group $H$ with associated subgroups $A$ and $B$. If the Search
Membership  Problem is decidable for subgroups $A$ and $B$ in $H$
then Algorithm I finds a reduced form for every given $g\in G$.
\end{proposition}

\subsection{Algorithm II for computing  normal forms}

Let the Search Membership Problem  be  decidable in $H$ for the
subgroups $A$ and $B$. Assume  also that the Coset Representative
Search Problem ($\mathbf{ CRSP}$) is decidable for the subgroups $A$
and $B$ in $H$, that is, there exist recursive sets $S_A$ and $S_B$
of representatives of $A$ and $B$ in $H$ and two algorithms
 which for a given word $w \in F(X)$ find, correspondingly,  a representative for
$Aw$ in $S_A$ and for $Bw$ in $S_B$. Notice that if $s_w$ is the
representative of $Aw$ in $S_A$ then $ws_w^{-1} \in A$, so, applying
to $ws_w^{-1}$  the algorithm for the Search Membership Problem for
$A$, one can find a representation of $w$ in the form $w = as_w$,
where $a$ is an element of $A$ given as a product of the generators
of $A$.

 Now we describe
the standard Algorithm II for computing normal forms of elements in
$G$.

Algorithm II  can be viewed as a sequence of applications  of
rewriting rules of the type
 \begin{itemize}
  \item $t^{-1}h \rightarrow \phi(c)t^{-1}s, $ \ where $h = cs, \ c
  \in A, \ s \in S_A$;
   \item $th \rightarrow \phi^{-1}(c)ts, $ \  where $h = cs, \ c
  \in B, \ s \in S_B$;
  \item $t^\epsilon t^{-\epsilon} \rightarrow 1$
  \end{itemize}
to a given element $g \in G$ presented as a word in the standard
generators of $G$. Since the problems $\mathbf{ SMP}$ and
$\mathbf{CRSP}$ are decidable for $A$ and $B$ in $H$ the rewriting
rules above are effective (i.e., given the left side of the rule one
can effectively find the right side of the rule). The rewriting
process is organized
 ``from the right to the left'', i.e, the algorithm always rewrites {\em the rightmost
 occurrence}   of the left side of a rule above.

It is not hard to see that the Algorithm II halts on every input $g
\in G$ in finitely many steps and yields the normal form of $g$.

We summarize the discussion above in the following  theorem.

\begin{theorem}
\label{th:normal-forms}
 Let $G = \left< H, t \mid t^{-1}At = B
\right>$ be an HNN-extension of a group $H$ with associate subgroups
$A$ and $B$. If the Search Membership Problem  and  the Coset
Representative Search Problem are decidable for subgroups $A$ and
$B$ in $H$ (with respect to fixed transversals $S_A$ and $S_B$) then
Algorithm II
 finds the normal form  for every given $g \in G$.
\end{theorem}

If  elements in the group $H$ admit some particular "normal form",
then one can define a normal form for elements of $G$. Namely, let
$\nu(h)$ be a normal form  of an element $h \in H$ - usually we
assume that $\nu(h)$ is a particular word in the generators of $H$
 uniquely representing  the element $h$.  Then the $\nu$-normal form of an
 element $g = h_0t^{\epsilon_1}h_1 \cdots t^{\epsilon_n}h_n \in G$ given
 as in (\ref{eq:normalforms}) is defined by
 $$
\nu(g) = \nu(h_0)t^{\epsilon_1}\nu(h_1) \cdots
t^{\epsilon_n}\nu(h_n).
 $$
 In this case  $\nu(g)$ is a word in the generators of $G$  uniquely
 representing the element $g$. Notice, that if there is an algorithm
 to compute $\nu$-forms of elements in $H$ then there is an
 algorithm to compute $\nu$-forms of elements in $G$ provided $G$
 satisfies the conditions of Theorem \ref{th:normal-forms}.

\subsection{Algorithm III for computing cyclically reduced normal
forms} \label{sec:Algorithm_II}

Now we want to briefly outline an algorithm which, given an element
$g \in G$ in reduced form, computes its cyclically reduced normal
form. Recall that the cyclically reduced normal form  of $g$ is a
conjugate of $g$ whose  normal form is cyclically reduced.  We work
under the assumption that the Search Membership Problem  and the
Coset Representative Search Problem  are decidable for
 subgroups $A$ and $B$ in
$H$, so one can use the standard Algorithm II to find normal forms
of elements of $G$. Assume now that the  {\em Conjugacy Membership
Search Problem} ($\mathbf{CMSP}$) is  also decidable for subgroups
$A$ and $B$ in $H$. The latter means that for a given $g \in H$ one
can determine whether or not $g$ is a conjugate of an element from
$A$ (or from $B$), and if so, find such an element in $A$ (in $B$)
and a conjugator.

\bigskip
\noindent {\sc Algorithm III: Computing Cyclically Reduced Normal
Forms.}

\medskip \noindent {\sc Input:} a word in the reduced form
\[
g = h_0 t^{\epsilon_1} h_1 \cdots h_{k-1}t^{\epsilon_k} h_k.
\]

\begin{itemize}
\item[{\sc Step 0}] Find the normal form of $g$ using Algorithm II:
$$
g = hp_1\cdots p_k.
$$

 \item[{\sc Step 1}]\
 \begin{itemize}
\item If $l(g) = 0$ then $g \in H$.
 \begin{itemize}
  \item If $g \in C$, where where
$C=A\cup B$, or if $g$ is not
 conjugate to an element in $C$, then $g$ is already in
cyclically reduced form.
 \item If $g^x \in C$ for some $x \in H$ then use a decision
 algorithm for \textbf{CMSP} to find a particular such $x$  and
 replace $g$ by $g^x$.
\end{itemize}

\item If $l(g) = 1$, then $g$ is already in cyclically reduced
form.

\item If $l(g) \geqslant 2$ and $\epsilon_1 = \epsilon_k$ then $g$
is already in cyclically reduced form.

\end{itemize}

\item[{\sc Step 2}]\

\textsc{If} $l(g) \geqslant 2$ and $\epsilon_1 = -\epsilon_k$ and
$s_kh \not\in A$ (when $\epsilon_k =-1$) or $t_kh \not\in B$ (when
$\epsilon_k =1$) then $g$ is in cyclically reduced form.

\textsc{Otherwise}, if $s_kh \in A$ then set \[ g^*
=t^{-\epsilon_1}h^{-1}ght^{\epsilon_1}.\] Obviously, we have $l(g^*)
= l(g)-2$,  and we can input $g^*$ to Step 0 and iterate.

The case $t_kh \in B$ is treated similarly.

\end{itemize}

\begin{theorem}
\label{th:cyc-red-forms} Let $G = \left< H, t \mid t^{-1}At = B
\right>$ be an HNN-extension of a group $H$ with associate subgroups
$A$ and $B$. If the Search Membership Problem,   the Coset
Representative Search Problem, and the   Conjugacy Membership Search
Problem  are decidable for subgroups $A$ and $B$ in $H$  then
Algorithm III
 finds the cyclically reduced normal form  for every given $g \in G$.
\end{theorem}

\section{The Conjugacy search problem for regular elements}

\subsection{The Conjugacy criterion}
\label{sec:conjugacy-criterion}

In this section we formulate, in a slightly modified form, the
well known  conjugacy criterion for HNN-extensions, due to Collins
\cite{collins}.

Observe, that any element of $G$ has a conjugate of the type $h_0
t^{\epsilon_1} \cdots h_{r-1}t^{\epsilon_r}$.  Recall that  the
\emph{$i$-cyclic permutation} of a cyclically reduced element
 $$g =
h_0 t^{\epsilon_1} \cdots h_{r-1}t^{\epsilon_r}
 $$
is the element
\[
g_i =h_i t^{\epsilon_{i+1}} \cdots t^{\epsilon_{r}} h_0
t^{\epsilon_{1}} \cdots h_{i-1} t^{\epsilon_{i}},
\]
rewritten in normal form.

\begin{theorem}
\label{th:Collins} Let $G = \left< H, t \mid t^{-1}At = B \right>$
be an HNN-extension of the base group $H$ with associated
subgroups $A$ and $B$. Let
\[
g = h_0 t^{\epsilon_1} \cdots h_{r-1}t^{\epsilon_r}, \quad g' =
h'_0 t^{\eta_1} \cdots h'_{s-1}t^{\eta_s}
\]
 be conjugate cyclically reduced elements of $G$. Then one of the
following is true:
\begin{itemize}

\item [1)] Both $g$ and $g'$ lie in the base group $H$. If $g \not\in A
\cup B$ then $g' \not\in A \cup B$ and $g$ and $g'$ are conjugate
in $H$.

\item [2)] If $g\in A\cup B$ then $g' \in A\cup B$ and there exists a
finite sequence of elements $c_1,\dots, c_l \in A\cup B$, such
that $c_0 = g$, $c_l = g'$ and $c_i$ is conjugated to $c_{i+1}$ by
an element of the form $ht^\epsilon$, $h\in H$, $\epsilon = \pm
1$.

\item [3)] Neither of $g$, $g'$ lies in the base group $H$, in which
case $r=s$ and $g'$ can be obtained from $g$ by $i$-cyclically
permuting it\/ {\rm ($i = 1,\dots, r$)} and then conjugating it by
an element $z$ from $A$, if\/ $\epsilon_i = -1$, or from $B$, if\/
$\epsilon_i = +1$.
\end{itemize}

\end{theorem}

\subsection{Bad pairs}

 Let $C = A\cup B$. We say that $(c,g) \in C \times G$ is a \emph{bad}
pair if $c \ne 1$, $g\not\in C$, and $gcg^{-1} \in C$. We will show
later that bad pairs is the main source of "hardness" of the
Conjugacy  Problem in $G$.

The following lemma gives a more detailed description of bad pairs.

\begin{lemma}
Let $c \in C \smallsetminus \{1\}$,   $g \in G \smallsetminus C$,
and let $g = hp_1\cdots p_k$ be the  normal form of $g$. Then
 $(c,g)$ is a bad pair if
and only if the following system of equations has solutions $c_1,
\dots, c_{k+1}\in C$.
\begin{eqnarray*}
p_kcp_k^{-1} &=& c_1\\
p_{k-1}c_1p_{k-1}^{-1} &=& c_2\\
&\vdots & \\
p_{1}c_{k-1}p_{1}^{-1} &=& c_k\\
hc_kh^{-1} &=& c_{k+1}.
\end{eqnarray*}
 \label{lm:equation}
\end{lemma}

\begin{proof} This lemma is a special case of Lemma~\ref{lm:gg'} below.
\end{proof}

\medskip

We denote the system of equations in Lemma~\ref{lm:equation} by
$B_{c,g}$. Observe that the consistency of the system $B_{c,g}$ does
not depend on the particular choice of representatives of $A$ and
$B$ in $H$. Sometimes we shall treat $c$ as a variable, in which
case the system will be denoted $B_g$. If $c,c_1, \ldots, c_{k+1}
\in C \smallsetminus \{1\} $ is a solution of $B_g$  then we call it
a {\em nontrivial} solution of $B_g$.

 Now we want to study slightly more general equations of
the type $gc =c'g'$ and their solutions $c,c'\in C$.

\begin{lemma} Let $G = \left< H,t \mid t^{-1}At = B \right>$. Let $g,g'\in
G$ be elements given by their normal forms \begin{equation} g =
hp_1\cdots p_k, \qquad g' = h'p_1'\cdots p_k'
\end{equation}
 Then the equation $gc = c'g'$ has a solution
$c,c' \in C$ if and only if the following system $S_{g,g'}$ of
equations in variables $c,c_1,\dots, c_k$ has a solution in $C$.
\begin{eqnarray*}
p_kc &=& c_1p_k'\\
p_{k-1}c_1 &=& c_2p_{k-1}'\\
&\vdots& \\
p_1c_{k-1} &=& c_kp_1'\\
hc_k &=& c'h'
\end{eqnarray*}
\label{lm:gg'}
\end{lemma}
 \begin{proof}
The proof of Lemma~\ref{lm:gg'} is a word-by-word reproduction of
the proof of Lemma~4.5 in \cite{amalgam-1}.
 \end{proof}
The first $k$ equations of the system $S_{g,g'}$ form what we call
the \emph{principal system of equations}, we denote it by
$PS_{g,g'}$. In what follows we consider $PS_{g,g'}$ as a system
in variables $c,c_1,\dots, c_k$ which take values in $C$, the
elements $p_1,\dots, p_k, p'_1,\dots,p'_k$ are constants.

\subsection{Regular elements and Black holes}

The set

\[
N^*_G(C) = \left\{ g \mid C^g \cap C \ne 1 \right\}
\]
is called  the {\em generalized normalizer} of the set $C$.

Notice that if $(c,g)$ is a bad pair then $g \in N^*_G(C)
\smallsetminus C$ and $c \in Z_g(C)$, where
\[
Z_g(C) = \{\, c \in C \mid c^{g^{-1}} \in C\,\} = C^g \cap C.
\]
 We refer to the  set
\[
\BH = N^*_G(C)
\]
as to a \emph{Black Hole} of the Conjugacy Problem in $G$. Elements
from $\BH$ are called \emph{singular}, and elements from $R = G
\smallsetminus \BH$ \emph{regular}. The following description of the
black hole is an immediate corollary of Lemma~\ref{lm:equation}.

\begin{corollary}
Let $G = \left< H, t \mid t^{-1}At = B \right>$. Then an element
 $g \in G \smallsetminus C$ is singular if and only if the system
$B_{g}$ has a nontrivial solution $c,c_1,\dots,c_{k+1} \in C$.
\end{corollary}

\begin{lemma}
\label{le:criterion-regularity}
 Let\/ $G = \left< H,t \mid t^{-1}At
= B \right>$ and $g,g' \in G$. If\/ $l(g)=l(g') \geqslant 1$ and the
system $PS_{g,g'}$ has more than one solution in $C$ then the
elements $g,g'$ are singular.
\end{lemma}
 \begin{proof}
 The proof repeats the proof of Lemma~4.10 in \cite{amalgam-1}.
 \end{proof}

\subsection{Effective recognition of regular elements}

Let $M$ be a subset of a group $G$. If $u,v \in G$, we call the set
$uMv$ a \emph{$G$-shift} of $M$. For a collection ${\mathcal{M}}$ of
subsets in $G$, we denote by $SI({\mathcal{M}},G)$ the least set of
subsets of $G$ which contains $\mathcal{M}$ and is closed under
$G$-shifts and finite intersections.

\begin{lemma}
Let $G$ be a group and $C = A\cup B$ be the union of two subgroups
$A$ and $B$ of $G$. If $D \in SI(\{C\},G)$ and $D \ne \emptyset$
then $D$ is the union of finitely many sets of the form
\[ D = (A^{g_1}\cap \cdots \cap A^{g_m} \cap B^{g'_1}\cap \cdots
\cap B^{g'_n})h
\]
for some elements $g_1,\dots,g_m,g'_1,\dots,g'_n,h \in G$.
\label{lem:3}
\end{lemma}

\begin{proof}
The proof of this lemma repeats the proof of Lemma~4.7 of
\cite{amalgam-1}. \end{proof}

\begin{corollary}
Let $Sub(C)$ be the set of all finitely generated subgroups of $C$.
Then every  non-empty set from $SI(Sub(C), H)$ is a finite union of
cosets of the type
\begin{equation} \label{eq:SI}
(C^{h_1}\cap \cdots \cap C_m^{h_m})h,
\end{equation}
where $C_i \in Sub(C)$ and $h_i, h \in H$.
\end{corollary}
 Now we can apply these results to solution sets of the systems $PS(g,g')$.
\begin{lemma}
\label{le:E}
 Let\/ $G = \left< H,t \mid t^{-1}At = B \right>$. Then
for any two elements $g$ and\/ $g'$ with normal forms
\[
g = hp_1\cdots p_k, \quad g' = h'p'_1\cdots p'_k \quad (k
\geqslant 1)
\]
the set\/ $E_{g,g'}$ of all elements $c$ from $C$ for which the
system $PS_{g,g'}$ has a solution $c,c_1,\dots,c_k$, is equal to
\[
E_{g,g'}  = C \cap p_k^{-1}Cp'_k \cap \cdots \cap p_k^{-1}\cdots
p_1^{-1} C p'_1\cdots p'_k.
\]
In particular, if\/ $E_{g,g'} \ne \emptyset$ then it is the union of
at most\/ $2^{k+1}$ cosets of the type $A^{g_1}\cap \cdots \cap
A^{g_m} \cap B^{g'_1}\cap \cdots \cap B^{g'_n}$, where
$g_1,\dots,g_m,g'_1,\dots,g'_n \in G$.
\end{lemma}
 \begin{proof}
The proof of this lemma is essentially the same as that of
Lemma~4.8 in \cite{amalgam-1}.
 \end{proof}

 Following
\cite{amalgam-1} we say that the {\em Cardinality Search Problem} is
decidable for  $SI(Sub(C),H)$ if for a set $D \in SI(Sub(C),H)$,
given as a finite union of sets obtained by finite sequences of
shifts and intersections from subgroups from $Sub(C)$,  one can
effectively decide whether $D$ is empty, finite, or infinite and, if
$D$ is finite non-empty, list all elements of $D$. It is not hard to
see that the Cardinality Search Problem for $SI(Sub(C),H)$ is
decidable if and only if one can decide the Cardinality Problem for
intersections of cosets of the type $C_1h_1 \cap C_2h_2$, where
$C_1,C_2 \in Sub(C), h_1, h_2 \in H$.

\begin{corollary}
\label{co:eff-E}
 Let $G = \left< H,t \mid t^{-1}At = B \right>$. If
the Cardinality Search Problem is decidable in $SI(Sub(C),H)$, then,
given $g,g'$ as above, one can effectively find the set $E_{g,g'}$.
In particular, one can effectively check whether  $E_{g,g'}$ is
empty, singleton, or infinite.
\end{corollary}
\begin{proof}The proof repeats the proof of Corollary~4.9 in
\cite{amalgam-1}. \end{proof}

\begin{theorem}
Let $G = \left< H,t \mid t^{-1}At = B \right>$ be an HNN-extension
of a  finitely presented group $H$ with finitely generated
associated subgroups $A$ and $B$. Set\/ $C = A \cup B$. Assume also
that $H$ allows algorithms for solving the following problems:
 \bi
 \item The Search Membership Problem for $A$ and $B$ in $H$.
  \item The Coset Representative Search Problem for
subgroups $A$ and $B$ in $H$.
 \item The Cardinality Search Problem for $SI(Sub(C),H)$ in $H$.
\item The Membership Problem for $N^\ast_H(C)$ in $H$.
  \ei
    Then there exists
an algorithm for deciding whether a given element in $G$ is
regular or not. \label{lm:regular}
\end{theorem}

\begin{proof} For a given $g \in G$ we can find the canonical
normal  form of $g$ using Algorithm II. There are two cases to
consider:

1) If $l(g) \geq 1$ then by Lemma \ref{le:criterion-regularity} $g$
is singular if and only if the system $B_{c,g}$ has a nontrivial
solution $c, c_1, \ldots, c_k \in C$. Now, if $B_{c,g}$ has no
solutions in $C$ (and we can check it effectively) then $g$ is
regular. If $B_{c,g}$ has precisely one solution then we can find it
and check whether it is trivial or not, hence we can find out
whether $g$ is regular or not. If $B_{c,g}$ has more then one
solution (and we can verify this effectively) then $g$ is not
regular, since if the system $B_{c,g}$ has two distinct solutions
then one of them is nontrivial.

2) If $l(g) = 0$ then $g \in H$. In this case $g$  is regular if and
only if $g \not \in N^\ast_H(C)$.  Since the Membership problem for
$N^\ast_H(C)$ is decidable in $H$ one can check if $g$ is regular or
not. \end{proof}

\begin{corollary}
Let $G = \left< H,t \mid t^{-1}At = B \right>$ be an HNN-extension
of a free group  $H$ with finitely generated associated subgroups
$A$ and $B$. Then the set of regular elements in $G$ is recursive.
\end{corollary}

\subsection{The Conjugacy search problem and regular elements}

The aim of this section is to study the Conjugacy Search Problem for
regular elements in HNN-extensions. Recall that the Conjugacy Search
Problem is decidable in a group $G$ if there exists an algorithm
that for given two elements $g,h \in G$ decides whether these
elements are conjugate in $G$ or not, and if they are the algorithm
finds a conjugator. We show that the conjugacy search problem for
regular elements is solvable under some very natural restrictions on
the group $H$. We start with the following particular case of the
Conjugacy Search Problem.

\begin{theorem}
\label{th:conjugacy-reg-non-zero-length}
 Let $G = \left< H,t \mid t^{-1}At = B
\right>$ be an HNN-extension of a finitely presented group $H$ with
associated finitely generated subgroups $A$ and $B$.
 Assume also that $H$
allows algorithms for solving the following problems:
 \bi
 \item The Word Problem in $H$.
 \item The Search Membership Problem for $A$ and $B$ in $H$.
 \item The Coset Representative Search Problem for subgroups $A$ and $B$
in $H$.
\item The Cardinality Search Problem for $SI(Sub(C),H)$ in $H$.
\ei Then the Conjugacy Search Problem in $G$ is decidable for
arbitrary pairs $(g,u)$, where $g$ is an element that has a
cyclically reduced regular normal form of non-zero length, and $u
\in G$.
\end{theorem}
 \begin{proof}
 Let $g$ be a fixed regular cyclically reduced element of length $l(g) \geqslant 1$
 and $g'$ be an arbitrary element from $G$.  By Theorem
 \ref{th:normal-forms} one can find the normal form of $g$
 and this normal form is cyclically reduced. If $l(g') = 0$ then
 by the Conjugacy criterion $g'$ is not a conjugate of $g$. Suppose now that $l(g') \geq 1$.
 The decidability of the Search  Membership Problem and
the Coset Representative Search Problem for subgroups $A$ and $B$ in
$H$ allows one to apply  Steps 1 and 2 of Algorithm  III to $g'$
until either the element $g^\ast$
 in the Step 2 becomes of length $0$, or the element $g^\ast$  becomes a cyclically
 reduced normal form of length $l(g^\ast) \geq 1$.
 In the former case $g'$ is not a conjugate of $g$, and in the later case we found a
 cyclically  reduced normal form of $g'$ with  $l(g') \geq 1$.
 For simplicity, we may assume now  that the elements  $g$ and $g'$
 are in cyclically reduced normal forms:
$$ g = cp_1 \ldots p_k, \ \ \   g^{\prime} = c^{\prime}p_1^{\prime}
 \ldots  {p_{k^{\prime}}}^{\prime}.$$
 According to the Conjugacy criterion, the elements $g$
and $g'$ are conjugate in $G$ if and only if $k = k^\prime$ and for
some $i$-cyclic permutation $\pi_i(g')$ of  $g'$ the equation
$c^{-1}gc = \pi_i(g')$ has a solution $c$ in $C$.  By Lemma
\ref{lm:gg'} the equation $c^{-1}gc = \pi(g')$ has a solution in $C$
if and only if the system $S_{g, \pi(g')}$ has a solution in $C$.
Since $g$ is regular the system $PS_{g, \pi(g')}$ has at most one
solution in $C$. Decidability of the Cardinality Search Problem
problems for $SI(Sub(C),H)$ in $H$ allows one to check whether
$PS_{g, \pi(g')}$ has a solution in $C$ or not, and if it does, one
can find the solution (see Lemma \ref{le:E} and Corollary
\ref{co:eff-E}). Now one can verify whether this solution satisfies
the last equation of the system $S_{g, \pi(g')}$ or not. If not, the
system   $S_{g, \pi(g')}$ has no solutions in $C$, as well as the
equation $c^{-1}gc = \pi(g')$. Otherwise, the system $S_{g,
\pi(g')}$  and the equation $c^{-1}gc = \pi(g')$ have solutions in
$C$ and we have found one of these solutions.

Suppose now that $g$ is an arbitrary element from $G$ that have  a
cyclically reduced regular normal form $g_1$ of length $l(g_1)
\geqslant 1$. We claim that one can find a conjugate $g_2$ of $g$
that has a  regular cyclically reduced non-zero normal form. Indeed,
using Algorithm III one can find a cyclically reduced normal form,
say $g'$ of $g$. Clearly, $l(g') \geq 1$. Observe that  $g'$ and
$g_1$ are conjugated in $G$, hence  by the Conjugacy criterion $g_1
= \pi_i(g')^c$ for some $i$-cyclic permutation $\pi_i(g')$ of $g'$
and some $c \in C$. Since $g_1$ is regular this implies that
$\pi_i(g') =  g_1^{c^{-1}}$ is also regular (easy calculation). It
follows that  one of cyclic permutations of $g'$ is regular. Now one
can effectively list all cyclic permutations $\pi_j(g')$ of $g'$ and
apply the decision algorithm described above to each pair
$\pi_j(g'),u)$. This proves the theorem.

 \end{proof}
Now we study the Conjugacy Search Problem for regular elements of
length $0$.

\begin{lemma}
Let $G = \left< H,t \mid t^{-1}At = B \right>$.    If the Search
Membership Problem
 and the Coset Representative Search Problem
for subgroups $A$ and $B$ in $H$ are decidable and the Conjugacy
Search Problem for $H$ is  decidable then the Conjugacy Search
Problem is decidable for all pairs of elements $(g,u)$ where $g$ is
a regular element  of $G$ with $l(g) = 0$ and $u$ is an arbitrary
element from $G$.
\end{lemma}
  \begin{proof}
Let $g \in G$ be a regular element of length $0$. It follows that $g
\in H$ and $g \not \in N^\ast_H(C)$. In particular, $g \in H
\smallsetminus (A \cup B)$.  As was mentioned in the proof of
Theorem \ref{th:conjugacy-reg-non-zero-length} the decidability of
the Search Membership Problem and the Coset Representative Search
Problem for subgroups $A$ and $B$ in $H$ allows one  either to check
that the cyclically reduced normal form of $u$ has length $\geq 1$
or to find a conjugate $u'$ of $u$ that belongs to $H$ as a word in
generators of $H$. In the former case, by the Conjugacy criterion
$u$ is not a conjugate of $g$. In the latter case,  again by the
Conjugacy criterion $g$ and $u$ are conjugate in $G$ if and only if
$g$ and $u'$ are conjugate in $H$. Using the decision algorithm for
the Search Conjugacy Problem in $H$ one can check if $g$ and $u'$
are conjugate in $H$, and if they are, find a conjugator. This
finishes the proof of the lemma.
  \end{proof}

\begin{corollary}
Let $G = \left< H,t \mid t^{-1}At = B \right>$ be an HNN-extension
of a free group $H$ with associated finitely generated subgroups $A$
and $B$. Then the Conjugacy Search Problem in $G$ is decidable for
arbitrary pair $(g,u)$ where $g$ is a regular element in $G$ and $u
\in G$.
\end{corollary}

\section{Miller's construction}
\label{se:Miller}

In this section we discuss  a particular type of HNN-extension
introduced by C.~Miller III in \cite{miller2}.

Let
$$H = \langle s_1, \ldots, s_n \mid R_1, \ldots, R_m \rangle$$
 be a finitely presented group. Starting with $H$ one can
 construct a new group $G(H)$, the Miller group of $H$,  with generators:
 \begin{equation} \label{eq:1-1}
 q,\  s_1, \ldots, s_n,\  t_1, \ldots, t_m, \ d_1, \ldots , d_n
 \end{equation}
and relators:
  \begin{equation}
  \label{eq:1-2}
  t_i^{-1}qt_i = qR_i, \ \ t_i^{-1}s_jt_i = s_j, \ \ d_j^{-1}qd_j =
 s_j^{-1}qs_j, \ \ d_k^{-1}s_jd_k = s_j
 \end{equation}
Generators from (\ref{eq:1-1}) are  called the {\em standard}
generators of $G(H)$.

One can realize $G(H)$ as
 a generalized mapping torus of a free group, which is  a very particular
  type of a multiple HNN-extension of a free group. To this end
put
 $$S = \{s_1, \ldots, s_n\}, \ \ D = \{d_1, \ldots, d_n\}, \ \
 T = \{t_1, \ldots, t_m\}$$
 and denote by $q$ a new symbol not in $S \cup T \cup D$. Let
  $$F(S,q) = F(q, s_1, \dots, s_n)$$
be a free group with basis $S_q = \{q\} \cup S$.

For every $i = 1, \ldots, m$ we define an automorphism $\phi_i$ of
$F(S,q)$ as
 \[
  \phi_i: \left \{
                \begin{array}{rcl}
                q & \rightarrow & qR_i \\
                s_j & \rightarrow & s_j
                \end{array}
         \right.
     \]
For every $k = 1, \ldots, n$ we define an automorphism $\psi_k$ of
$F(S,q)$ as
\[
  \psi_k: \left \{
                \begin{array}{rcl}
                q & \rightarrow & s_k^{-1}qs_k \\
                s_j & \rightarrow & s_j
                \end{array}
         \right.
     \]
 It is easy
to see that  the following multiple  HNN-extension of $F(S,q)$ with
the stable letters from $T \cup D$ has precisely the same
presentation   (\ref{eq:1-2}) as the group $G(H)$ in the standard
generators, so it is isomorphic to $G(H)$:
 \begin{equation}
 \label{eq:1-3}
G(H) \simeq   \langle F(S,q),   T \cup D  \mid \
 t_i^{-1} f t_i = \phi_i(f), \ \ d_k^{-1}fd_k = \psi_k(f), \ \ f \in F(S,q)\ \rangle
  \end{equation}

 As noticed in \cite{miller1}, the group $G(H)$ can be also viewed as
the standard HNN-extension of a direct product of two free groups by
a single stable letter $q$. Indeed, consider the following
construction.

The subgroup
$$\langle T \cup D\rangle  \leqslant G(H)$$
is free with a basis $T \cup D$ (since its image in the quotient
group of $G(H)$ modulo  the normal closure of $F(S,q)$ is free), we
denote it by $F(T,D)$. The subgroup $\langle S \rangle$ of $G(H)$ is
also free with basis $S$ (as a subgroup of $F(S,q)$, which, in its
turn,  is a subgroup of $G(H)$), we denote it by $F(S)$.

Put
 $$K = F(T,D) \times F(S).$$
Then the following are free subgroups of $K$:
\begin{eqnarray*}
 A &=& \langle t_1, \ldots, t_m, s_1d_1^{-1}, \ldots,
 s_nd_n^{-1}\rangle,\\
 B &=& \langle t_1R_1^{-1}, \ldots, t_mR_m^{-1}, s_1d_1^{-1}, \ldots,
 s_nd_n^{-1}\rangle.
 \end{eqnarray*}
They are isomorphic under the map

\[
  \theta: \left \{
                \begin{array}{lcll}
                t_i & \rightarrow & t_iR_i^{-1}, & i = 1, \ldots, m \\
                s_jd_j^{-1} & \rightarrow & s_jd_j^{-1}, & i = 1, \ldots,
                n
                \end{array}
         \right.
     \]

It is a straightforward verification that the following
HNN-extension of $K$ with the stable letter $q$ and the subgroups
$A, B$ associated via $\theta$ has precisely the same presentation
(\ref{eq:1-2}) as the group $G(H)$ in the standard generators, so
it is isomorphic to $G(H)$:
 \begin{equation}
 \label{eq:1-4}
 G(H) \simeq \langle K, q \mid q^{-1}aq = \theta(a)  \mbox{ for }  a \in A \
 \rangle .
\end{equation}

Below we collect some elementary properties of $G(H)$.

\begin{lemma}
\label{le:semidirect} In this notation,
\begin{itemize}
 \item[{\rm (i)}] $\langle S \cup \{q\} \rangle  \simeq F(S,q)$,

 $\langle T \cup D \cup S \rangle  \simeq
 K$;
 \item[{\rm (ii)}] $F(S,q)$ is normal in $G(H)$;
 \item[{\rm (iii)}] $K = A \ltimes F(S)$, where $\ltimes$, as usual, denotes the semidirect product;
 \item[{\rm (iv)}] $K = B \ltimes F(S)$.
 \end{itemize}
 \end{lemma}

\begin{proof} Straightforward verification. \end{proof}

\begin{corollary}
\label{co:transversal} The set $F(S)$ is a system of left {\rm
(}and right\/{\rm )} representatives of $K$ modulo $A$, as well as
modulo $B$.
\end{corollary}

It follows from the definition of $K$ and Lemma
\ref{le:semidirect} that every element $x \in K$ can be uniquely
written in three different forms:
 \begin{equation}
 \label{eq:decomp1-2-3} x = u(x)s(x) = a(x)s_a(x) =
b(x)s_b(x),
 \end{equation}
  where $s(x), s_a(x), s_b(x) \in F(S)$, $u(x) \in F(T,D)$, $a(x) \in A$, $b(x) \in B$.

 \medskip
 \textsc{Convention:} All the groups that appeared above came equipped
 with particular sets of generators. From now on we fix these
 generating sets and call them {\em standard} generating sets.
 Furthermore, for all algorithms that we discuss below we assume
 that all  elements of  our groups, when these elements are viewed as inputs of the algorithms,
 they are presented as words in
 the standard generators  or their inverses. The same assumption is required
 for outputs of the algorithms.  Moreover, in this event we denote
 by $|g|_L$ the length of the word which represents $g$ in the standard generators of a group $L$.
 Instead of $|g|_{G(H)}$  we write $|g|$.

\begin{lemma}
\label{le:formssemidirect} For a given $x \in K$ one can
effectively find all three decompositions
 $$x = u(x)s(x) = a(x)s_a(x) = b(x)s_b(x),$$
in time at most quadratic in $|x|$. Moreover, the following
equalities hold for some constant $c$:
 \begin{enumerate}
 \item[{\rm (i)}] $|u(x)| \leqslant |x|$, $|s(x)| \leqslant |x|, $
 \item[{\rm (ii)}] $|a(x)|_A \leqslant |x|, \ \ |s_a(x)| \leqslant c \cdot |x|^2$,
 \item[{\rm (iii)}] $|b(x)|_B \leqslant |x|, \ \ |s_b(x)| \leqslant c\cdot |x|^2$,
 \end{enumerate}
\end{lemma}

\begin{proof} Let $x \in K$. To decompose $x$ into the form $x =
u(x)s(x)$ one needs only to collect in $x$ all letters from $(T
\cup D)^{\pm 1}$ to the left and all letters from $S^{\pm 1}$ to
the right.

To decompose $x$ in the form  $x= a(x)s_a(x)$ one can replace each
occurrence of the symbol $d_i^{-1}$ by $s_i^{-1}(s_id_i^{-1})$ and
each occurrence of $d_i$ by $(d_is_i^{-1})s_i$. This allows one to
present $x$ as a word in the standard generators of $A$ and
$F(S)$. Now, using the standard procedure for semidirect products
(and the relations from (\ref{eq:1-2}))   one can collect the
generators of $A$ to the left, which yields the result. Similar
 argument  provides an algorithm to present $x$ in the form
 $x = b(x)s(x)$. \end{proof}

\begin{corollary}
\label{co:unique-2} In the notations above the following hold:
 \begin{enumerate}
 \item[{\rm (i)}] For every $u \in F(T,D)$ there exists a unique $s \in F(S)$
such that $us \in A$. Moreover, one can find such $s$ in quadratic
time of $|u|$.
 \item[{\rm (ii)}] For any $g, h \in K$, if $u(g) = u(h)$ then $a(g) =
 a(h)$ and $b(g) = b(h)$.
 \end{enumerate}
\end{corollary}

\begin{proof} (i) comes directly from Lemmas \ref{le:semidirect} and
\ref{le:formssemidirect}. Now  (ii) follows from (i).
 \end{proof}

\subsection{Normal forms of elements of $G(H)$}

In this section we discuss normal  and cyclically reduced normal
forms  of elements of $G(H)$. We start with the standard normal
forms in HNN-extensions and then simplify  them using specific
properties of $G(H)$.

In what follows we view the group $G(H)$ as an HNN-extension of
the group $K$ by a single stable letter $q$:
 $$G(H) = \langle K, q \mid \ q^{-1}aq = \theta(a) \ \ (a \in
 A)\ \rangle $$
By Corollary  \ref{co:transversal} we  can choose the set $F(S)$
as the set of representatives of $K$ modulo $A$, as well as modulo
$B$. The general theory of HNN-extensions tells one that in this
event every element $g \in G(H)$ can be uniquely written in the
form
 \begin{equation}
 \label{eq:normalstandard}
 g = hq^{\varepsilon_1} s_1 \cdots q^{\varepsilon_k} s_k,
 \end{equation}
  where $s_i \in F(S)$, $\varepsilon_i \in \{1, -1\}$, $h
  \in K$, $ k \geqslant 0$, and if $\varepsilon_{i+1} = -\varepsilon_i$ then $s_i \neq 1$.
  Since $K = F(T\cup D) \times
  F(S)$ we can write $h$ uniquely  as a product $h = us_0$ where $u
  \in F(T,D)$ and $s_0 \in F(S)$. It follows that $g$ can be
  written uniquely as
   \begin{equation}
   \label{eq:newform}
   g = us_0q^{\varepsilon_1} s_1 \cdots q^{\varepsilon_k} s_k.
    \end{equation}
 We refer to (\ref{eq:newform})  as
to the \emph{normal form} of $g$.  Taking in account that $$f =
s_0q^{\varepsilon_1} s_1 \cdots q^{\varepsilon_k} s_k \in F(S,q)$$
one can rewrite (\ref{eq:newform}) in the form
\begin{equation}
 g = uf, \ \mbox{  where } \  u
\in F(T,D) \ \mbox{ and }\ f \in F(S,q). \label{eq:uf}
  \end{equation}

\begin{lemma}
\label{le:complexity} Let $H$ be a finitely presented group and
$G(H)$ be the corresponding Miller group. Then the following
conditions hold:
 \begin{enumerate}
  \item[{\rm (i)}] There is an algorithm  which for every element $g
  \in G(H)$  finds its
  normal form $(\ref{eq:newform})$. Moreover it
   has at most cubic time complexity in the length
  $|g|$.
 \item[{\rm (ii)}] Algorithm III {\rm (}which finds, for every element $g
  \in G(H)$, a  cyclically reduced
   element $g^\prime \in G(H)$ which is a conjugate of $g${\rm )},
   has at most cubic time complexity in the length
  $|g|$.
\end{enumerate}
\end{lemma}

\begin{proof} To prove  (i) we show that a slight   modification of the
standard Algorithm II
 does the job.  Let
  $$g = w_1q^{\epsilon_1}w_2q^{\epsilon_2} \cdots
  q^{\epsilon_n}w_{n+1},$$
   where $w_i \in K$, $\epsilon_i \in \{1,-1\}$.
Assume (by induction on $n$)  that one can effectively rewrite, in
at most $C_1\cdot 2n \cdot |v|^2$ steps, the word
 $$v = w_2q^{\epsilon_2} \cdots  q^{\epsilon_n}w_{n+1}$$
into the normal form
 $$v = u_2s_2q^{\epsilon_2} s_3 \cdots  s_nq^{\epsilon_n}s_{n+1}$$
where $u_2 \in F(T,D)$, $s_i \in F(S)$ and such that
 $$|u_2| \leqslant |v|, \ \ \ |s_i| \leqslant C_2|v|^2$$
 for some constant $C_2$ independent of $g$. Then
  $$g =w_1q^{\epsilon_1} v = w_1q^{\epsilon_1}u_2s_2q^{\epsilon_2} s_3 \cdots  s_nq^{\epsilon_n}s_{n+1}$$
Suppose, for certainty, that $\epsilon_1 = -1$ (the case
$\epsilon_1 = 1$ is similar). Now by Lemma
  \ref{le:formssemidirect} one can effectively rewrite $u_2$ in
  the form $as_a$ with $a \in A, s_a \in F(S)$ such that
   $$|a|_A \leqslant |u_2| \leqslant |v|, \ \ |s_a| \leqslant c|u_2|^2 \leqslant c|v|^2$$
   where $c$ is the constant from Lemma \ref{le:formssemidirect}.
   This rewriting requires at most $C_3|u_2|^2$ steps, where $C_3$
   is a constant from Lemma \ref{le:formssemidirect} which is independent of $u_2$.
 Then $q^{-1}a = \theta(a)q^{-1}$ and  $|\theta(a)|_B = |a|_A  \leqslant |v|$.
    Observe that
     $$|\theta(a)| \leqslant C_R|\theta(a)|_B \leqslant C_R|v|,$$
      where $C_R = \max \{\,|R_i| \ \mid \ i = 1, \dots, m\,\}.$
    Hence $|w_1\theta(a)| \leqslant |w_1|+ C_R|v| \leqslant C_R|g|$.
Again by Lemma
  \ref{le:formssemidirect} one can effectively rewrite $w_1\theta(a)$ in
  the form $us_1$ (in at most $C_R^2|g|^2$  steps)  where $u \in
  T(D,T)$, $s_1 \in F(S)$ and
   $$|s_1| \leqslant C_R^2|g|^2.$$
To estimate the length of $u$ notice that $u =
u(w_1)u(\theta(a))$, so
 $$|u| \leqslant |u(w_1)| + |u(\theta(a))|.$$
Observe that
 $$|u(\theta(a))| \leqslant |\theta(a)|_B  = |\theta(a)|_A \leqslant |v|.$$
 Hence
  $$|u| \leqslant |u(w_1)| + |u(\theta(a))| \leqslant |w_1| + |v| \leqslant |g|,$$
as required.  This argument shows how to find the normal form of
$g$ in the case when $q^{\epsilon_1}s_2q^{\epsilon_2}$ is not a
pinch. In the case when it is  a pinch one needs also to cancel
$q^{\epsilon_1}q^{\epsilon_2}$. In both cases the required bounds
on the length of elements are satisfied.  The total number of
steps required to write down the normal form of $g$ is bounded
from above by
 $$ C_1\cdot 2n \cdot |v|^2 + C_3|u_2|^2 +  C_R^2|g|^2 $$
 If we assume that $C_1 \geqslant  C_3$, $C_R$ then one can continue the chain of inequalities:
  $$\leqslant C_1 (2n|v|^2 + |v|^2 +|g|^2) \leqslant C_1 \cdot 2(n+1) \cdot |g|^2, $$
as required.

(ii) follows easily from (i) if $g \not \in K$. If $g \in K$ then
one has to verify whether $g \in A \cup B$ or not, and if yes,
then find a conjugate element in $A \cup B$. Since $K$ is a direct
product of two free groups the problem above reduces to the
Conjugacy Membership Problem \cite{amalgam-1} for finitely
generated subgroups of free groups which is decidable in at most
quadratic time (see \cite{km}).  This proves the lemma.
 \end{proof}

\subsection{Regular elements in $G(H)$} \label{sec:6}

In this section we show that even though the standard black hole
$\BH$ of $G(H)$ (given as an HNN-extension of $K$)  is very big, in
fact, it is equal to the whole group $G(H)$, one still can  show
that just a relatively small portion of elements of $\BH$  are
``hard'' for the conjugacy problem in $G(H)$. We refer to such
elements as to \emph{strongly singular}. On the contrary, the
elements for which the conjugacy problem is relatively easy are
called \emph{weakly regular}; see precise definitions below.

The following result shows that the standard black hole in $G(H)$
with respect to two different presentations of $G(H)$ as an
HNN-extension is  the whole group, and, as a result, the standard
notion of a regular element becomes vacuous.

\begin{lemma}
\label{le:BH} Let $G(H)$ be the Miller group of $H$. Then the
following hold:
\begin{itemize}
\item[{\rm (a)}] Let $G(H)$ be presented as the HNN-extension {\rm
(\ref{eq:1-3})} then
$$\BH = G(H).$$

\item[{\rm (b)}] Let $G(H)$ be presented as the HNN-extension {\rm
(\ref{eq:1-4})} of the group $K$ with the stable letter $q$ then
 $$\BH = G(H).$$
\end{itemize}
\end{lemma}

\begin{proof} Set $C= A\cup B$. It immediately follows from presentations
(\ref{eq:1-3}) and (\ref{eq:1-4}) and Lemma~\ref{le:semidirect} that
in the both cases $N^*_G(C) = G$. Since $\BH = N^*_G(C)$ the lemma
follows immediately. \end{proof}

Therefore we have to weaken the definition of regular elements. A
cyclically reduced element $g \in G(H)$ is called {\em weakly
regular} if in its normal form (\ref{eq:uf}) the element $u$ in
the decomposition $g = uf$ is non-trivial. If $u = 1$ then $g$ is
called \emph{strongly singular}.

We define the {\em strong black hole}  $\SBH(G)$  of $G(H)$ as the
set of all elements conjugate to strongly singular elements,
\[
\SBH(G) = \bigcup_{g\in G(H)} F(S,q)^g=F(S,q),
\]
for $F(S,q)$ is a normal subgroup in $G(H)$. Observe that every
cyclically reduced element in $G \smallsetminus \SBH(G)$ is weakly
regular.

 The main result of this section  is the following theorem.
\begin{theorem}
\label{th:conjugacy-criterion-regular} Let
$$g = uf = us_0q^{\epsilon_1} \cdots s_kq^{\epsilon_k}$$
be a weakly regular cyclically reduced element of\/ $G(H)$ and
$g^\prime = u^\prime f^\prime$ be an arbitrary  cyclically reduced
element of $G(H)$. If
 $$g^x = g^\prime$$
   for some $x = vh \in G(H)$ with $v \in F(T \cup D)$ and $h \in F(S,q)$  then the following conditions hold:
\begin{enumerate}
 \item[{\rm (i)}] $g^\prime$ is weakly regular and $u^v = u^\prime$. Therefore, replacing
 $g^\prime$ by $(g^\prime)^{v^{-1}}$  and\/ $x$ by $xv^{-1}$ we may assume
 that\/ $u^\prime = u$ and\/ $x = h \in F(S,q)$.
 \item[{\rm (ii)}] If  $g \in K \smallsetminus (A \cup B)$ {\rm (}that is, $f \in F(S)${\rm )} then
$f^\prime \in F(S)$ and $f^s = f^\prime$ for some $s \in F(S)$.

 \item[{\rm (iii)}] If\/ $g \in A \cup B$ then  $g^\prime \in
A \cup B$. Moreover, the following hold:
 \begin{enumerate}
 \item[{\rm (iii.a)}] If\/ $g$ and\/ $g^\prime$ are in
the same factor then $g = g^\prime$.\/
  \item[{\rm (iii.b)}] If\/ $g \in A$ and $g^\prime
\in B$ then $q^{-1}gq = g^\prime$.
 \item[{\rm (iii.c)}] If\/ $g \in B$ and\/ $g^\prime \in
A$ then $qgq^{-1} = g^\prime$.
\end{enumerate}
 \item[{\rm (iv)}] If\/  $g \not \in K$ then $g^\prime \not \in K$ and there
 exists an $i$-cyclic permutation
$$g^\ast = us_0^\prime q^{\epsilon_1} \cdots s_k^\prime q^{\epsilon_k}$$
of $g^{\prime}$ and an element $z \in A \cup B$ such that
$$g^z = g^\ast ,$$
and $z \in A $ if\/ $\epsilon_k = -1$, and\/ $z \in B$ if\/
$\epsilon_k = 1$. Moreover, in this case there exists an integer\/
$l$ and elements $y, c \in F(S)$  such that:
\begin{enumerate}
\item[{\rm (iv.a)}] $z= u_0^ly^l$ \ \ where $u_0$ is a generator
of the cyclic centralizer $C(u)$ in the group $F(D \cup T)$;
 \item[{\rm (iv.b)}] $$q^{-1}u_0yq = u_0c,  \ \ \ if \ \epsilon_k = -1$$
$$qu_0yq^{-1} = u_0c,  \ \ \ if \ \epsilon_k = 1$$
\item [{\rm (iv.c)}] If\/ $k = 1$ then
\begin{equation}
 \label{eq:55}
 s_{0}^\prime = y^{-l} s_{0} c^l,
\end{equation}

\item[{\rm (iv.d)}] If\/  $\epsilon_{k-1}\epsilon_k = 1$  then
\begin{equation}
 \label{eq:6}
 s_{k}^\prime = y^{-l} s_{k} c^l,
\end{equation}
If\/ $\epsilon_{k-1}\epsilon_k = -1$  then
\begin{equation}
 \label{eq:7}
 s_{k}^\prime = c^{-l} s_{k} c^l,
\end{equation}

\end{enumerate}
\end{enumerate}
\end{theorem}

\begin{proof} (i) Since $F(S,q)$ is normal in $G(H)$ (Lemma
\ref{le:semidirect}) one has
$$ u^\prime f^\prime = g^\prime = g^x = (uf)^{vh} = (u^vf^v)^h = u^v(f^v[u^vf^v,h]) $$
 where $u^v \in F(T,D)$ and $f^v[u^vf^v,h]  \in F(S,q)$. Uniqueness of the normal forms implies
$u^\prime = u^v$ and $f^\prime = f^v[u^vf^v,h]$. The equality $g^x =
g^\prime$ implies $g^{xv^{-1}} = (g^\prime)^{v^{-1}}$ hence
replacing $x$ by $xv^{-1} = vhv^{-1} \in F(S,q)$  and $g^\prime$ by
$(g^\prime)^{v^{-1}}$ one can assume that $g^\prime = uf^\prime$ and
$x \in F(S,q)$. This proves (i).

(ii)  follows immediately from the first case of the Conjugacy
Criterion (Theorem~\ref{th:Collins} in
Section~\ref{sec:conjugacy-criterion}) and from the decomposition
of $K$ into a direct sum of free groups
$$K= F(D \cup T) \times F(S)$$

(iii) Recall that every  element $g \in K$ can be decomposed
uniquely as $g = u(g)s(g)$ where $u(g) \in F(T,D)$, $s(g) \in F(S)$
(see Section \ref{se:Miller}). Now let $g \in A \cup B$. In this
event by the Conjugacy Criterion $g ^\prime \in A \cup B$. Since $x
\in F(S,q)$ then (as was shown above)
 $$u(g) = u(g^x) = u(g^\prime).$$
 By Corollary \ref{co:unique-2} this  implies
  $$a(g) = a(g^\prime), \ \ b(g) = b(g^\prime).$$
Therefore if $g$ and $g^\prime$ are in the same factor then  $g =
g^\prime$. If $g \in A$ and $g^\prime \in B$ then $q^{-1}gq =
g^\prime$. Indeed, in this case $g = a(g)$ and $a(g)^q = b(g) =
b(g') = g'$ since $g' \in B$.  Similarly, if $g \in B$ and $g^\prime
\in A$ then $qgq^{-1} = g^\prime$. This proves (iii).

(iv) By the Conjugacy Criterion if  $g \not \in K$ then $g^\prime
\not \in K$ and there
 exists an $i$-cyclic permutation
$$g^\ast = us_0^\prime q^{\epsilon_1} \cdots s_k^\prime q^{\epsilon_k}$$
of $g^{\prime}$ and an element $z$ such that
$$g^z = g^\ast.$$
Furthermore, in this case  $z \in A $ if $\epsilon_k = -1$, and $z
\in B$ if $\epsilon_k = 1$. This proves the first part of (iv).

By the argument in (i) $z = u_1  s$ where $[u,u_1] = 1$ and $s \in
F(S)$. Observe that the group  $F(D \cup T)$ is free, and $u \neq 1$
(since $g$ is weakly regular) therefore $C(u) = \langle u_0 \rangle$
for some  $u_0 \in F(D \cup T)$ which is not a proper power. Hence
$u_1 = u_0^l$ for some $l \in \mathbb{Z}$. Replacing $u_0$ by
$u_0^{-1}$ we may assume that $l
> 0$. It follows from Lemma \ref{le:semidirect} that $s= y^l$ for
some uniquely defined $y \in F(S)$. So $z = u_0^ly^l$ and  (iv.a)
follows.

The equality  $g^z = g^\ast$ implies $gz = zg^\ast$ which amounts to
 \begin{equation}
 \label{eq:gxxgprime}
  us_0q^{\epsilon_1} \cdots s_kq^{\epsilon_k} u_0^l y^l = u_0^l y^l
  us_0^\prime q^{\epsilon_1} \cdots s_k^\prime q^{\epsilon_k}.
  \end{equation}
   If $\epsilon_{k} = -1$ then  there exists $c \in F(S)$
  such that
   $$q^{-1}u_0yq = u_0c.$$
Similarly, if $\epsilon_{k} = 1$ then  there exists $c \in F(S)$
  such that
   $$qu_0yq^{-1} = u_0c.$$
This shows (iv.b).

 Rewriting now the left hand side of
(\ref{eq:gxxgprime}) into  normal form and comparing to the right
hand side of (\ref{eq:gxxgprime}) one can see that the following
equalities hold in the free group $F(S)$:

If $k = 1$ then:
 $$ s_{0}^\prime = y^{-l} s_{0} c^l,$$
and the case (iv.c) follows.

If  $k \geqslant 2$ then we have two subcases.

If $\epsilon_k = -1$ and $\epsilon_{k-1} = -1$, or if $\epsilon_k
= 1$ and $\epsilon_{k-1} = 1$ then:

 \begin{equation}
 s_{k}^\prime = y^{-l} s_{k} c^l,
\end{equation}

If $\epsilon_k = -1$ and $\epsilon_{k-1} = 1$, or if $\epsilon_k =
1$ and $\epsilon_{k-1} = -1$ then:

\begin{equation}
 s_{k}^\prime = c^{-l} s_{k} c^l,
\end{equation}

 This proves  (iv.d), and finishes the proof of the theorem.
\end{proof}

 \subsection{Conjugacy search problem in $G(H)$}

 The following result connects the
conjugacy problem in $G(H)$ with the word problem in $H$.

\begin{theorem}[Miller \cite{miller1}]
If the word problem is undecidable in $H$ then the Conjugacy Problem
is undecidable in $G(H)$.
\end{theorem}

This result shows that for strongly singular elements in $G(H)$ even
the classical decision form of the Conjugacy Problem is undecidable.
It turns out, however,  that  for weakly regular elements even the
Search Conjugacy Problem is decidable in $G(H)$. This  result
completes the general algorithmic picture of the Conjugacy Problem
in $G(H)$, even though one could still show that for many strongly
singular elements the Search Conjugacy Problem is decidable. We
leave for the future a more detailed analysis of the black hole
$\BH$ of $G(H)$.

  \begin{theorem}
  \label{th:main}
  Let $H$ be a finitely presented group and $G(H)$ be Miller's group based on $H$.
  Then the Conjugacy Search  Problem for pairs $(g,u)$, where $g$ is a weakly
  regular element from $G(H)$ and $u$ is an element from $G(H)$, is decidable in cubic time.
  \end{theorem}

\begin{proof} Let $g \in G(H)$ be a weakly regular element  of $G(H)$ and
$g^\prime$ be an arbitrary element of $G(H)$.

By Lemma \ref{le:complexity}, Algorithm III provides
us with the 
normal cyclically reduced forms $g = uf$ and
$g^\prime = u^\prime f^\prime$  in at most cubic time
 in the lengths $|g|$ and
$|g^\prime|$.

 In the rest of the proof  starting with cyclically reduced forms of elements
$g$ and $g^\prime$ we algorithmically verify whether or not the
cases (i)-(iv) of the Conjugacy Criterion
(Theorem~\ref{th:conjugacy-criterion-regular}) hold for these
elements. Simultaneously, we estimate time complexity of the
algorithm.

\paragraph{Case (i)}  One can easily check (in quadratic   time on $|u|
+ |u^\prime|$) whether or not the elements $u$ and $u^\prime$ are
conjugate in the free group $F(T,D)$. Moreover, if they are
conjugate then one can effectively find (in quadratic time on $|u|
+ |u^\prime|$)  a conjugator $v$.

Now we need to show that one can effectively write down the
element $(g^\prime)^{v^{-1}}$ in the normal form. Clearly,
 it suffices
to show on how one can effectively rewrite $(f^\prime)^{v^{-1}}$ as
a reduced word from $F(S,q)$.

 Using relations
$$q^{d_i} = s_i^{-1}qs_i,\quad  q^{t_i} = qR_i,\quad s_j^{t_i} = s_j,\quad s_j^{d_i} = s_j.$$
from the presentation  (\ref{eq:1-3}) of $G(H)$  one can rewrite
$(f^\prime)^{v^{-1}}$ as a word of length at most
$|f^\prime||v|\max\{|R_i| \mid i = 1, \dots, m\}$ in generators
$S\cup\{q\}$, and then freely reduce it.

 This shows that one can effectively check whether or not  the case (i) of the
 Conjugacy Criterion holds for $g$ and $g^\prime$. Moreover, if it
  holds then one can effectively find a required  element $v$
 and then effectively replace $g^\prime$  by
 $(g^\prime)^{v^{-1}}$.

\paragraph{Case (ii)} To determine effectively whether Case (ii) holds or
not one needs, firstly, to check whether $g \in A \cup B$ or not.
This amounts to the Membership Problem  for finitely generated
subgroups in free groups, which is linear. Secondly, one has to
solve the conjugacy problem in a free group, which is decidable
and at most quadratic.

\paragraph{Case (iii)} This case is obvious in view of the Case ~(ii).

\paragraph{Case (iv)} Verification of Case (iv) splits into two subcases:
firstly, one needs to find effectively the elements $u_0, y$, and
$c$, and, secondly, one has to find the number $l$, or prove that
such $l$ does not exist.

Since the element $u \in F(D \cup T)$  is  given, it is easy to
find its maximal root $u_0 \in F(D \cup T)$ in  quadratic time in
$|u|$.  Then by Lemma \ref{le:formssemidirect} one can find the
unique $y$ such that $u_0y \in A$ or $u_0y \in B$ (depending on
the sign of $\epsilon_k$). It takes again at most quadratic time.

Now one can effectively find the element  $c$ to satisfy (iv.b).
It follows again from Lemma \ref{le:formssemidirect}.

It is left to show  how one can effectively solve the systems in
(iv.c) in the free group $F(S)$ for an unknown $l$.

More generally, consider the following  equation in a free group
$F(S)$
$$a^lb^l = d$$
where $a, b,d \in F(S)$ are given, and $l$ is unknown  integer $l$.
In the {\em degenerate} case, where $d= 1$ and $a = b^{-1}$,  every
integer $l$ is a solution. Otherwise, this equation has at most one
solution in $F(S)$. Indeed,  if
 $$a^lb^l = d = a^mb^m$$
Then $a^{m-l} = b^{l-m}$ and $m = l$.

 Now we show how one can find
this unique solution if it exists. Below for elements $x,y,z \in
F(S)$  we write $x = y \circ z$ if $|x| = |y| + |z|$, i.e., no
cancellation in $yz$.

 If $[a,b] = 1$ then the equation takes the form $(ab)^l = d$ which is easy.
  Let  $[a,b] \neq 1$. We may assume that $a$ is cyclically reduced  and
$b= e^{-1} \circ b_0 \circ e$ for some $e, b_0 \in F(S)$ with $b_0$
cyclically reduced (one can find such $e, b_0$ in quadratic time).
There are three cases to consider.

If $ab = a \circ b$ then
 $$a^lb^l = a^l\circ e^{-1} \circ b_0^l \circ e = d$$
 hence
  $$l = \frac{|d|-2|e|}{|a| + |b_0|}.$$
If $e^{-1}$ does not cancel completely in $a^l e^{-1}$ then
  $a = a_1\circ a_2$, $e^{-1} = a^p \circ a_2^{-1}\circ e_1^{-1}$
   for some $a_1,a_2, e_1$. In this case
   $$a^lb^l = a^{l-p}\circ a_1 \circ e_1^{-1} \circ b_0^l \circ e_1
   \circ a_2 \circ a^{-p} = d$$
   and comparing length one can compute $l$ as before (since the elements
   $a_1, a_2, e_1$ are unique and can be easily found).

If $e^{-1}$ cancels completely in $a^l e_1^{-1}$ then the key fact
is that for any integers $k, m$  the cancellation in $a^k b_0^m$
cannot be longer than $|a| + |b_0|$ (otherwise the elements $a$
and $b_0$ (hence $a$ and $b$) commute). Again, one can make an
equation as above and solve it for $l$. We omit details here.

The argument above shows that one can find all possible values for
$l$ and then check whether the equation $a^lb^l = d$ holds in
$F(S)$. This requires at most quadratic number of steps.

Now, if  the elements $g$ and $g^\prime$ fall into premises of one
of the cases (ii) or (iii) then they are conjugate in $G(H)$ if and
only if the corresponding conditions, stated in the cases (ii) and
(iii),  hold. In this case the conjugator $x$ is easy to find.

If the elements $g$ and $g^\prime$ fall into premises of the cases
(i) or (and)  (iv), but  the corresponding conditions, stated in
these cases, do not hold in $G(H)$, then $g$ and $g^\prime$ are not
conjugate in $G(H)$.

If  $g$ and $g^\prime$ fall into the premises of the cases (i) and
(iv) and the cases  hold in $G(H)$ then one can effectively find the
unique solution $l$ of the systems in (iv.c), (iv.d) provided the
system is non-degenerate (defined above and described in more
details below). Hence the conjugating element $z$ (if it exists)
must be equal to $u_0^ly^l$. Now using the normal form algorithm one
can check whether, indeed, $g^z = g^\ast$ for $z = u_0^ly^l$ or not.

Finally, suppose that the equations (\ref{eq:55}), (\ref{eq:6}), and
(\ref{eq:7}) are degenerate. Observe, that  the equations
(\ref{eq:55}) and (\ref{eq:6}) can be written as
 $$(y^{-1})^l(s_ics_i^{-1})^l = s_i^\prime s_i^{-1}, \ \ \ i = 0,
 k$$
 in which case  they are degenerate if and only if $s_i^\prime = s_i$ and $y =
 s_ics_i^{-1}$.  For the equation  (\ref{eq:55}) (Case  (iv.c))  this implies that  $g = g^\ast$
 and $z = 1$. For the equation (\ref{eq:6}) (Case (iv.d), if
 $\epsilon_{k-1}\epsilon_k = 1$)  the following equalities hold in
 the event of $\epsilon_k = -1$ (the case $\epsilon_k = 1$ is
 similar and we omit it):
\begin{eqnarray*}
 g^z &=& (us_0s_{k-1}q^{\epsilon_{k-1}})^z(s_kq^{-1})^z\\ &=& (us_0s_{k-1}q^{\epsilon_{k-1}})^z
  u_0^{-l}y^{-l}
  s_ku_0^lc^lq^{-1}\\ &=&(us_0s_{k-1}q^{\epsilon_{k-1}})^z
  s_kq^{-1}.
\end{eqnarray*}
 Hence $g^z = g^\ast$ is equivalent to
  $$(us_0s_{k-1}q^{\epsilon_{k-1}})^z = u_0s_0^\prime
  q^{\epsilon_1} \cdots s_{k-1}^\prime q^{\epsilon_{k-1}}.$$
   This allows one to find $z$ by induction on $k$.

In the case of (\ref{eq:7}) (Case (iv.d), if
$\epsilon_{k-1}\epsilon_k = -1$) one has $s_k^\prime = s_k$ and
$c^{-1}s_k c = s_k$. Hence (in the case of $\epsilon_k = -1$)
  $$qs_kq^{-1}z = qs_kq^{-1}u_0^ly^l = qu_0^ls_kc^lq^{-1} =
  q\theta(z)c^{-l}s_kc^lq^{-1} = zqs_kq^{-1}.$$
Now $g^z = g^\ast$ is equivalent to
 $$(us_0q^{\epsilon_1} \cdots s_{k-1})^z = us_0^\prime
 q^{\epsilon_1} \cdots s_{k-1}^\prime$$
and, again, one can find $z$ by induction on $k$.

 This completes the proof of the theorem. \end{proof}

\section{Some algorithmic and probabilistic estimates}

\subsection{Asymptotic density}
\label{subesc:measuring}

 In this section we use the terminology and techniques
developed in \cite{BMS,multiplicative,KMSS,KMSS2} for measuring
various subsets of a  free group $F$ of rank $n$. This gives the
asymptotic classification of the sizes of these subsets.

Let $R$ be a subset of the free group $F$ and \[S_k = \{\, w \in F
\mid |w|=k\,\}
\]
the sphere of radius $k$ in $F$. The fraction
\[f_k(R)= \frac{\vert R \cap S_k\vert}{\vert S_k\vert}
\]
is the relative frequency of elements from $R$ among the words of
length $k$.  $R$ is called {\em generic}  if
 $\rho(R) = 1$, where  the {\em asymptotic density}
 $\rho(R)$ is defined by
  $$\rho(R) = \limsup_{k \rightarrow \infty} f_k(R).$$
 If, in addition, there
exists a positive constant $\delta <1$ such that $1-\delta^k <
f_k(R) < 1$ for all sufficiently large $k$ then $R$ is called
\emph{strongly generic}.

A set $R \subseteq F$ is {\em negligible} (\emph{strongly
negligible})  if its complement $F - R$ is generic (strongly
generic).

\subsection{The Strong Black Hole in Miller groups}

The Strong Black Hole $\SBH(G)$ in $G(H)$ is the subgroup $F(S,q)$
(see Section \ref{sec:6}).

\begin{theorem} \label{thm:8.1}
Let $H=\left<s_1,\ldots, s_n \mid R_1, \ldots, R_m \right>$ be a
finitely presented group and $G(H)$ the Miller group of $H$. Assume
that $m > 1$. Then the Strong Black Hole $\SBH(G)$ in $G(H)$ is a
strongly negligible set  and
\[
f_k(\SBH(G))<\left(\frac{n+1}{n+m}\right)^{k-1}, \hbox{ for all }
k>1.
\]
\end{theorem}

Note that the restriction $m>1$ is natural in the context of this
paper since one relator groups have decidable word problem by the
classical result of Magnus.

\begin{proof} Denote by $G_k, B_k, P_k$ the sets of all elements of length
$k$ the groups $G$, $F(S,q)$ and $F(T,D)$ respectively. Then it
follows from Equation (\ref{eq:uf}) that if $g=uf$ with $u\in
F(T,D)$ and $f \in F(S,q)$ then $l(g)=l(u)+l(f)$. This implies that:
\[
|G_k|=|P_k|+|P_{k-1}||B_1|+\cdots+|B_k|.
\]
Consequently, for $m>1$, we have
 \bea f_k(\SBH(G)) & = &
\frac{|B_k|}{|G_k|} < \frac{|B_k|}{|P_k|}\\
 &=&
\frac{(2n+2)(2n+1)^{k-1}}{(2n+2m)(2n+2m-1)^{k-1}}\\
& < & \left(\frac{n+1}{n+m}\right)^{k-1}. \eea \end{proof}

Below we present a quantitative estimate for the group $G(H)$, in
the case when $H$ is a well-known group with undecidable word
problem.

\begin{example}
Borisov constructed a group {\rm (}see \cite{col2}{\rm )} with
undecidable word problem with $10$ generators and $27$ relations:
\begin{gather}
\begin{split}
G=& \langle a, b,c,d,e,p,q,r,t,k \mid p^{10}a=ap,\; p^{10}b=bp,\;
p^{10}c=cp,\; p^{10}d=dp,\\  p^{10}e&=ep,\;
 qa=aq^{10},\; qb=bq^{10},\;
qc=cq^{10},\; qd=dq^{10},\; qe=eq^{10},\; ra=ar,\\ rb&=br,\;
rc=cr,\; rd=dr,\;
 re=er,\;  pacqr=rpcaq,\; p^2adq^2r=rp^2daq^2,\\
 p^3&bcq^3r=rp^3cbq^3,\;
p^4bdq^4r=rp^4dbq^4,\; p^5ceq^5r=rq^5ecaq^5,\;
\\ p^6de&q^6r=rp^6edbq^6,\;
p^7cdcq^7r=p^7cdceq^7,\; p^8caaaq^8r=rp^8aaaq^8,\\
 p^9daaaq^9r&=rp^9aaaq^9,\; pt=tp,\; qt=tq,\; k(aaa)^{-1}t(aaa) =
(aaa)^{-1}t(aaa)k\rangle.
\end{split}
\end{gather}
In our case $$\frac{n+1}{n+m}=\frac{11}{37}<\frac{1}{3}.$$
 Then, for instance,
$$f_{81}(\SBH(G)) < \left(\frac{n+1}{n+m}\right)^{80}<{\frac{1}{3^{80}}},$$
a number small beyond any practical possibility to find an element
in $\SBH(G)$ by picking random elements in $G$.
\end{example}

\subsection{Random elements in the base group}

In view of the general conjugacy criterion for HNN-extensions
(Theorem~\ref{th:Collins}), the most challenging case of the
Conjugacy Problem for Miller group $G(H)$ given  in the form
(\ref{eq:1-4})
$$
 G(H) \simeq \langle K, q \mid q^{-1}aq = \theta(a)  \mbox{ for }  a \in A \
 \rangle
$$
  is presented by pairs $(g,g')$ where
both elements $g$ and $g'$ belong to the base group $K$.

Let us look at random elements in $K$ using the measure-theoretic
framework of \cite{multiplicative}. A natural way to introduce an
atomic measure on $K$ is to use the direct sum decomposition $K =
F(T,D) \times F(S)$ and set
\[
\mu(k) = \mu_{\sigma_1}(u)\mu_{\sigma_2}(s)
\]
where $k = (u,s)$, $u \in F(T,D)$ and $s \in F(S)$, and
$\mu_{\sigma_1}$ and $\mu_{\sigma_2}$ are multiplicative measures
with stopping probabilities $\sigma_1$ and $\sigma_2$ on groups
$F(T,D)$ and $F(S)$, correspondingly (see Appendix below).

\begin{theorem}
\[
P(k \mbox{ is strongly singular}) = \sigma_1,
\]
 where $\sigma_1$ is the stopping probability of the random
word generator for the group $F(T,D)$.
\end{theorem}

\begin{proof} Let $k = us$ with $u \in F(T,D)$ and $s \in F(S,q)$. Since
$\SBH(G) \cap K = F(S)$,  it follows immediately that the element
$k$ belongs to $\SBH(G)$ if and only if $u = 1$. Hence the
probability in question is the probability $P (u=1) = \sigma_1$.
\end{proof}

\subsection{Definition of a measure on $G(H)$}
\label{subsec:negligible}  Similarly, one can introduce a measure
$\mu$ on the whole  Miller group $G(H)$. Indeed,  any element $g$
from $G(H)$ can be written uniquely as follows
\begin{equation} \label{eq:*}
g=uf, \hbox{ where } u\in F(T,D) \hbox{ and } f\in F(S,q).
\end{equation}

Let $\mu_1$ and $\mu_2$ be atomic measures defined for free groups
$F(T,D)$ and $F(S,q)$. Then an atomic measure $\mu$ for $G(H)$ is
defined on  $g$ by
\[
\mu(g)=\mu_1(u) \mu_2(f).
\]

This measure is quite natural and allows one to estimate sizes of
various sets of elements in $G(H)$.

\section{Appendix: Measuring sets in free groups}

\subsection{Generation of random words}
\label{se:1-2}

For completeness of exposition, we reproduce here some definitions
from \cite{multiplicative}.

 Let $F = F(X)$ be a free group with basis $X = \{x_1, \ldots,
x_m\}$. We use, as our random word generator, the following
no-return random walk on the Cayley graph $C(F,X)$ of $F$ with
respect to the generating set $X$.  We start at the identity element
$1$ and either do nothing with probability $s\in (0,1]$ (and return
value $1$ as the output of our random word generator), or move to
one of the $2m$ adjacent vertices with equal probabilities
$(1-s)/2m$. If we are at a vertex $v \ne 1$, we either stop at $v$
with probability $s$ (and return the value $v$ as the output), or
move, with probability $\frac{1-s}{2m-1}$, to one of the $2m-1$
adjacent vertices lying away from $1$, thus producing a new freely
reduced word $vx_i^{\pm 1}$. Since the Cayley graph $(C(F,X)$ is a
tree and we never return to the word we have already visited, it is
easy to see that the probability $\mu_s(w)$ for our process to
terminate at a word $w$ is given by the formula
\begin{equation}
\mu_s(w) = \frac{s(1-s)^{|w|}}{2m\cdot (2m-1)^{|w|-1}} \quad \hbox{
for } w \ne 1
\end{equation}
and
\begin{equation}
\mu_s(1) =s.
\end{equation}
Observe that the set of all  words of length $k$ in $F$ forms the
sphere $S_k$ of radius $k$ in $C(F,X)$ of cardinality  $|S_k| =
2m(2m-1)^{k-1}$. Therefore the probability to stop at a word of
length $k$ is
\begin{equation}
P(|w| = k) = s(1-s)^k. \label{eq:geometric}
\end{equation}
Hence the {\em lengths} of words produced by our process are
distributed according to a geometric law. It is obvious now that the
same random word generator can be described in simpler terms: we
make random freely reduced words $w$ of random length $|w|$
distributed according to the geometric law (\ref{eq:geometric}) in
such way that words of the same length $k$ are produced with equal
probabilities.

    The mean length $L_s$ of words in $F$
distributed according to $\mu_s$ is equal to $$L_s = \sum_{w \in F}
|w|\mu_s(w) = s \sum_{k= 1}^{\infty}k(1-s)^{k-1} = \frac{1}{s} -1.
$$ Hence we have a family of probability distributions $\mu =
\{\mu_s\}$ with the stopping probability $s \in (0,1)$ as a
parameter, which is related to the  average length $L_s$ as
$$s = \frac{1}{L_s+1}.$$
By  $\mu(R)$ we denote the function
\begin{eqnarray*}
\mu(R): (0,1) & \rightarrow & {\mathbb{R}}\\
s & \mapsto & \mu_s(R)= \sum_{w\in R} \mu_s(w);
\end{eqnarray*}
we call it  {\it measure} of $R$ with respect to the family of
distributions $\mu$.

Denote by $n_k = n_k(R) = |R \cap S_k|$ the number of elements
 of length $k$ in $R$, and by $f_k = f_k(R)$ the relative
frequencies
$$f_k = \frac{|R \cap S_k|}{|S_k|}$$ of words of length $k$ in
$R$. Notice that $f_0 = 1 $ or $0$ depending on whether $R$ contains
$1$ or not. Recalculating $\mu_s(R)$ in terms of $s$, we immediately
come to the formula $$ \mu_s(R) = s\sum_{k=0}^\infty f_k(1-s)^k, $$
and the series on the right hand side is convergent for all $s \in
(0,1)$. Thus, for every subset $R \subseteq F$, $\mu(R)$ is an
analytic function of $s$.

The asymptotic behaviour of the set $R$ when $L_s \rightarrow
\infty$ depends on the behaviour of the function $\mu(R)$ when
$s\rightarrow 0^+$.   Here we just mention how one can obtain a
first coarse  approximation of the asymptotic behaviour of the
function $\mu(R)$.   Let $W_0$ be the  no-return non-stop random
walk  on $C(F,X)$ (like $W_s$ with $s = 0$), where the walker moves
from a given vertex  to any adjacent vertex away from the initial
vertex 1 with equal probabilities $1/2m$. In this event, the
probability $\lambda(w)$ that the walker hits an element $w \in F$
in $|w|$ steps (which is the same as the probability that the walker
ever hits $w$)  is equal to
$$
\lambda(w) = \frac{1}{2m(2m-1)^{|w|-1}}, \ \hbox{ if } \ w \neq 1, \
\ \hbox{ and }  \ \lambda(1) = 1.
$$
 This gives rise to an atomic measure
$$ \lambda(R) = \sum_{w \in R}\lambda(w) = \sum_{k=0}^{\infty}f_k(R)$$
where $\lambda(R)$  is just the sum of the relative frequencies of
$R$.
 This measure is not probabilistic,
since some sets have no finite measure (obviously, $\lambda(F) =
\infty$), moreover, the  measure $\lambda$ is finitely additive, but
not $\sigma$-additive. We shall call $\lambda$  the {\em frequency }
measure on $F$. If $R$ is $\lambda$-measurable (i.e., $\lambda(R) <
\infty$) then $f_k(R) \rightarrow 0$ when $k \rightarrow \infty$, so
intuitively, the set $R$ is ``small" in $F$.

 A number of papers (see, for example, \cite{Arzh,Fr,Olsh,woess}), used the {\em asymptotic density} (or more,
precisely, the {\em spherical asymptotic density})
 $$\rho(R) = \lim\sup f_k(R)$$
  as a numeric characteristic of the set
$R$ reflecting its asymptotic behavior.


 A more subtle analysis of asymptotic
behaviour of $R$ is provided  by  the
 {\em relative growth  rate} \cite{grigorchuk0}
$$
\gamma(R) = \lim\sup \sqrt[k]{f_k(R)}.
$$
Notice the obvious inequality $\gamma(R) \leqslant 1$. If $\gamma(R)
< 1$  then the series $\sum f_k$ converges. This shows that if
$\gamma(R) < 1$ then $R$ is $\lambda$-measurable.

\subsection{The multiplicativity of the measure and generating functions}
\label{se:1-3}

It is convenient to renormalise our measures $\mu_s$  and work with
the parametric  family $\mu^* = \{\mu_s^* \}$ of  {\em adjusted
measures}
\begin{equation}
\mu_s^*(w) = \left(\frac{2m}{2m-1} \cdot \frac{1}{s} \right)\cdot
\mu_s(w).
\end{equation}
This new measure $\mu_s^*$  is {\it multiplicative} in the sense
that
\begin{equation}
\mu_s^*(u\circ v) = \mu_s^*(u)\mu_s^*(v),
\end{equation}
where $u\circ v$ denotes the product of non-empty words $u$ and $v$
such that $|uv| = |u| +|v|$ i.e.\ there is no cancellation between
$u$ and $v$. The measure $\mu$ itself is {\em almost multiplicative}
in the sense that
\begin{equation}
\mu_s(u\circ v) = c\mu_s(u)\mu_s(v) \quad \hbox{for} \quad c =
\frac{2m}{2m-1} \cdot \frac{1}{s}
\end{equation}
for all non-empty words $u$ and $v$ such that $|uv| = |u| +|v|$.
Therefore our measure is close in its properties to the
\emph{Boltzmann samplers} of \cite{DFLS}: there, random
combinatorial objects are generated with probabilities obeying the
following rule: \emph{if thing $A$ is made of two things $B$ and $C$
then $p(A)= p(B)p(C)$}.

If we denote
\begin{equation}
t = \mu_s^*(x_i^{\pm 1})=  \frac{1-s}{2m-1} \label{eq:adjusted}
\end{equation}
 then
\begin{equation}
\mu_s^*(w) = t^{|w|}
\end{equation}
 for every non-empty word $w$.

Similarly, we can adjust the frequency  measure $\lambda$ making it
into a multiplicative atomic measure
\begin{equation}
\lambda^*(w) = \frac{1}{(2m-1)^{|w|}}. \label{eq:lambda*}
\end{equation}

Let now $R$ be a subset in $F$ and  $n_k = n_k(R) = |R \cap S_k|$ be
the number of elements of length $k$ in $R$.  The sequence
$\{n_k(R)\}_{k=0}^{\infty}$ is called the {\it spherical growth
sequence} of $R$. We assume, for the sake of minor technical
convenience, that $R$ does not contain the identity element $1$, so
that  $n_0 = 0$. It is easy to see now that $$ \mu^*_s(R) =
\sum_{k=0}^\infty n_kt^k .$$  One can view $\mu^*(R)$ as the
generating function of the spherical growth sequence of the set $R$
in variable $t$ which is convergent for each  $t \in [0,1)$. This
simple observation allows us  to apply a well established machinery
of generating functions of context-free languages \cite{flajolet} to
estimate probabilities  of sets.

\subsection{Cesaro density}
\label{se:1-4}

 Let $\mu = \{\mu_s\}$ be the parametric family of distributions defined
above. For a subset $R$ of  $F$ we define the {\it limit measure}
$\mu_0(R):$
 $$ \mu_0(R) = \lim_{s \rightarrow 0^+}\mu_s(R) =  \lim_{s \rightarrow 0^+} s \cdot \sum_{k=0}^\infty
f_k(1-s)^k.$$  The function $\mu_0$   is additive, but not
$\sigma$-additive, since $\mu_0(w) = 0$ for a single element $w$. It
is easy to construct a set $R$ such that $\lim_{s\rightarrow
0^+}\mu(R)$ does not exist. However, in the applications that we
have in mind we have not yet encountered such a situation. Strictly
speaking, $\mu_0$ is not a measure because the set  of all
$\mu_0$-measurable sets is not closed under intersections (though it
is closed under complements). Because $\mu_s(R)$ gives an
approximation of $\mu_0(R)$ when $s \rightarrow 0^+$, or
equivalently, when $L_s \rightarrow \infty$, we shall call $R$ {\it
measurable at infinity} if $\mu_0(R)$ exists, otherwise $R$ is
called {\it singular}.

If $\mu(R)$ can be expanded  as a convergent power series in $s$ at
$s=0$ (and hence in some neighborhood of $s = 0$):
 $$ \mu(R) = m_0 + m_1s + m_2s^2 + \cdots, $$ then
  $$\mu_0(R) = \lim_{s
\rightarrow 0^+} \mu_s(R) = m_0,$$ and an easy corollary from a
theorem by Hardy and Littlewood \cite[Theorem~94]{hardy}
asserts that $\mu_0$ can be computed as the {\em Cesaro limit}
 \begin{equation}
 \label{eq:cesaro-1}
  \mu_0(R)  = \lim_{n\rightarrow\infty} \frac{1}{n}\left(f_1+\cdots+f_n\right).
  \end{equation}
   So it will be also natural to call $\mu_0$ the {\em Cesaro density}, or
    {\it asymptotic average density}.

The Cesaro density $\mu_0$ is  more sensitive than the standard
asymptotic density $\rho= \lim\sup f_k$. For example,  if $R$ is a
coset of a subgroup $H$ of finite index in $F$ then it follows from
Woess \cite{woess} that $$ \mu_0(R) = \frac{1}{|G:H|}, $$ while,
obviously,  $\rho(H) = 1$ for the group $H$ of index $2$ consisting
of all elements of even length.

 On the other hand, if $\lim_{k\rightarrow \infty} f_k(R)$
exists (hence is equal to $\rho(R)$) then $\mu_0(R)$ also exists and
$\mu_0(R) = \rho(R)$. In particular, if a set $R$ is
$\lambda$-measurable, then it is $\mu_0$-measurable, and $\mu_0(R) =
0$.

\subsection{Asymptotic classification of subsets}
\label{se:1-5}

In this section we introduce  a classification of subsets $R$ in $F$
according to the asymptotic behaviour of the functions $\mu(R)$.

Let $\mu = \{\mu_s\}$ be the family of measures defined in Section
\ref{se:1-2}. We start with a  global characterization of subsets of
$F$.

Let $R$ be a subset of $F$. By its construction, the function
$\mu(R)$ is analytic on  $(0,1)$.
We say that $R$ is {\em smooth} if $\mu(R)$ can be analytically
extended to a neighborhood of $0$.

We start by considering a linear approximation of $\mu(R)$.
 If the set $R$ is smooth then  the linear term in the expansion  of
$\mu(R)$ gives a linear approximation of  $\mu(R)$:
$$
\mu_s(R) = m_0 + m_1s + O(s^2).
$$
 Notice that, in this case,  $m_0 = \mu_0(R)$ is the Cesaro density of $R$.
An easy corollary of \cite[Theorem~94]{hardy} shows
 that
 if $\mu_0(R) = 0$ then
$$ m_1 = \sum_{k=1}^\infty f_k(R) = \lambda(R). $$

On the other hand, even without assumption that $R$ is smooth,
 if $R$ is $\lambda$-measurable (that is, the series $\sum
f_k(R)$ converges), then
 $$\mu_0(R) = 0 \ \  \hbox{ and }\ \ \lim_{s\rightarrow 0^+} \frac{\mu(s)}{s} =  \lambda(R).$$
This allows us to use for the limit
$$\mu_1 = \lim_{s\rightarrow 0^+} \frac{\mu(s)}{s},$$
if it exists, the same term {\em frequency measure} as for
$\lambda$. The function $\mu_1$ is an
  additive   measure on $F$ (though it is not $\sigma$-additive).

 Now we can introduce a subtler classification of sets in $F$:
\begin{itemize}
\item {\em Thick} subsets: $\mu_0(R)$ exists, $\mu_0(R) > 0$ and
$$\mu(R) = \mu_0(R) + \alpha_0(s), \ \ where \ \ \lim_{s \rightarrow
0^+}\alpha_0(s) = 0.$$

\item {\em Sparse} subsets: $\mu_0(R) = 0$, $\mu_1(R)$ exists  and
$$\mu(R) = \mu_1(R)s + \alpha_1(s)\ \  where \ \ \lim_{s \rightarrow
0^+}\frac{\alpha_1(s)}{s}  = 0.$$

\item {\em Intermediate density} subsets: $\mu_0(R) = 0$
 but $\mu_1(R)$ does not exist.

\item {\em Singular\/} sets: $\mu_0(R)$ does not exist.

\end{itemize}

For sparse sets, the values of $\mu_1$ introduce a further and more
subtle discrimination by size.

It can be easily seen \cite{multiplicative} that every
$\lambda$-measurable set is sparse.

A set $R \subseteq F$ is \emph{strongly negligible} if there exist
positive $\delta <1$ such that $f_k(R) < \delta^k$. It is easy to
see that every strongly negligible set is sparse and
$\lambda$-measurable.

\bigskip

\normalsize

\vfill

 \noindent \textsf{Alexandre V. Borovik, School of Mathematics,  PO Box 88,
The University of Manchester, Sackville Street, Manchester M60 1QD,
United Kingdom}

\noindent {\tt borovik@manchester.ac.uk}

\noindent {\tt http://www.ma.umist.ac.uk/avb}

\medskip

\noindent \textsf{Alexei G. Myasnikov,
 Department of Mathematics and Statistics, McGill University, 805 Sherbrooke St W.,
 Montreal, QC, H3A 2K6, Canada}

\noindent  {\tt amiasnikov@gmail.com}

\noindent {\tt http://www.math.mcgill.ca/~alexeim/}

\medskip
\noindent \textsf{Vladimir N. Remeslennikov, Omsk Branch of
Mathematical Institute SB RAS,\linebreak 13 Pevtsova Street, Omsk
644099, Russia}

\noindent {\tt remesl@iitam.omsk.net.ru}

\end{document}